\documentclass[a4paper]{amsart}
\usepackage[all]{xy}
\usepackage{amssymb,enumitem,nicefrac,mathtools}
\usepackage[british]{babel}

\newcommand*{\MRref}[2]{\linebreak[0] \href{http://www.ams.org/mathscinet-getitem?mr=#1}{MR \textbf{#1}}}
\newcommand*{\arxiv}[1]{\linebreak[0] \href{http://www.arxiv.org/abs/#1}{arXiv:#1}}
\usepackage[pdftitle={C*-Algebras over Topological Spaces: the Bootstrap Class},
  pdfauthor={Ralf Meyer and Ryszard Nest},
  pdfsubject={Mathematics; MSC 19K35, 46L35, 46L80, 46M20}
]{hyperref}
\usepackage[lite]{amsrefs}
\usepackage{microtype}

\DeclareMathOperator{\Hom}{Hom}
\DeclareMathOperator{\coker}{coker}
\DeclareMathOperator{\Prim}{Prim}

\newcommand*{\KK}{\textup{KK}}

\newcommand*{\K}{\textup K}

\newcommand*{\congto}{\xrightarrow\cong}
\newcommand*{\into}{\rightarrowtail}
\newcommand*{\prto}{\twoheadrightarrow}

\newcommand*{\cl}[1]{\overline{#1}}% closure

\newcommand*{\primid}{\mathfrak p}% primitive ideal
\newcommand*{\Fil}{\mathfrak F}% filtration
\newcommand*{\CONT}{\textup C}% continuous functions
\newcommand*{\Hils}{\mathcal H}% Hilbert space or module
\newcommand*{\Ideal}{\mathcal I}% ideal in a triangulated category

\newcommand*{\Ab}{\mathfrak{Ab}} % category of Abelian groups
\newcommand*{\Cstarcat}{\mathfrak{C^*alg}}
\newcommand*{\Cstarsep}{\mathfrak{C^*sep}}
\newcommand*{\KKcat}{\mathfrak{KK}}

\newcommand*{\Cat}{\mathfrak C}  % generic category, often Abelian
\newcommand*{\pt}{\mathord\star}

\newcommand*{\bo}{\textup b}
\newcommand*{\maxi}{\textup{max}}
\newcommand*{\ID}{\textup{id}}

\newcommand*{\C}{\mathbb C}
\newcommand*{\R}{\mathbb R}
\newcommand*{\Z}{\mathbb Z}
\newcommand*{\N}{\mathbb N}
\newcommand*{\Comp}{\mathbb K}
\newcommand*{\Bound}{\mathbb B}

\newcommand*{\Ideals}{\mathbb I}% ideal lattice of a C*-algebra
\newcommand*{\Open}{\mathbb O}% open subsets of a top. space
\newcommand*{\Loclo}{\mathbb{LC}}% locally closed subsets of a top. space

\newcommand*{\nb}{\nobreakdash}% no break after this hyphen

\newcommand*{\loc}{\textup{loc}}
\newcommand*{\Bootstrap}{\mathcal B}
\newcommand*{\Bootstrapc}{\mathcal D}
\newcommand*{\Star}{$^*$\nb-}
\newcommand*{\Cst}{\textup C^*}

\newcommand*{\minotimes}{\otimes}
\newcommand*{\blank}{\text\textvisiblespace}
\newcommand*{\inOb}{\mathrel{\in\in}}

\newcommand*{\defeq}{\mathrel{\vcentcolon=}}

\newcommand*{\lad}{\vdash}

\numberwithin{equation}{section}

\theoremstyle{plain}
\newtheorem{theorem}[equation]{Theorem}
\newtheorem{proposition}[equation]{Proposition}
\newtheorem{lemma}[equation]{Lemma}
\newtheorem{corollary}[equation]{Corollary}
\theoremstyle{definition}
\newtheorem{definition}[equation]{Definition}
\theoremstyle{remark}
\newtheorem{example}[equation]{Example}

\begin{document}

\title[\(\Cst\)-Algebras over Topological Spaces: The Bootstrap Class]{\(\Cst\)-Algebras over Topological Spaces:\\ The Bootstrap Class}

\author{Ralf Meyer}
\address{Mathematisches Institut\\
         Georg-August-Universit\"at G\"ottingen\\
         Bunsenstra\ss e 3--5\\
         37073 G\"ottingen\\
         Germany\\
}
\email{rameyer@uni-math.gwdg.de}

\author{Ryszard Nest}
\address{K{\o}benhavns Universitets Institut for Matematiske Fag\\
         Universitetsparken 5\\ 2100 K{\o}benhavn\\ Denmark
}
\email{rnest@math.ku.dk}


\begin{abstract}
  We carefully define and study \(\Cst\)\nb-algebras over
  topological spaces, possibly non-Hausdorff, and review some
  relevant results from point-set topology along the way.  We
  explain the triangulated category structure on the bivariant
  Kasparov theory over a topological space and study the
  analogue of the bootstrap class for \(\Cst\)\nb-algebras over a
  finite topological space.
\end{abstract}

\maketitle

\section{Introduction}
\label{sec:intro}

If~\(X\) is a locally compact Hausdorff space, then there are
various equivalent characterisations of what it means for~\(X\)
to act on a \(\Cst\)\nb-algebra~\(A\).  The most common definition
uses an essential \Star{}homomorphism from \(\CONT_0(X)\) to
the centre of the multiplier algebra of~\(A\).  An action of
this kind is equivalent to a continuous map from the primitive
ideal space \(\Prim(A)\) of~\(A\) to~\(X\).  This makes sense
in general: A \(\Cst\)\nb-algebra over a topological space~\(X\),
which may be non-Hausdorff, is a pair \((A,\psi)\), where~\(A\)
is a \(\Cst\)\nb-algebra and \(\psi\colon \Prim(A)\to X\) is a
continuous map.  One of the purposes of this article is to
discuss this definition and relate it to other notions due to
Eberhard Kirchberg and Alexander Bonkat
\cites{Kirchberg:Michael, Bonkat:Thesis}.

An analogue of Kasparov theory for \(\Cst\)\nb-algebras over
locally compact Hausdorff spaces was defined already by Gennadi
Kasparov in~\cite{Kasparov:Novikov}.  He used it in his proof
of the Novikov conjecture for subgroups of Lie groups.
Kasparov's definition was extended by Eberhard Kirchberg to the
non-Hausdorff case in~\cite{Kirchberg:Michael}, in order to
generalise classification results for simple, purely infinite,
nuclear \(\Cst\)\nb-algebras to the non-simple case.  In his
thesis~\cite{Bonkat:Thesis}, Alexander Bonkat studies an even
more general theory and extends the basic results of Kasparov
theory to this setting.

This article is part of an ongoing project to compute the
Kasparov groups \(\KK_*(X;A,B)\) for a topological space~\(X\)
and \(\Cst\)\nb-algebras \(A\) and~\(B\) over~\(X\).  The aim is a
Universal Coefficient Theorem in this context that is useful
for the classification programme.  At the moment, we can
achieve this goal for \emph{some} finite topological spaces
(see~\cite{Meyer-Nest:Filtrated_K}), but the general situation,
even in the finite case, is unclear.

Here we describe an analogue of the bootstrap class for
\(\Cst\)\nb-algebras over a topological space.  Although we also
propose a definition for infinite spaces
in~\S\ref{sec:def_bootstrap_infinite}, most of our results are
limited to finite spaces.

\smallskip

Our first task is to clarify the definition of \(\Cst\)\nb-algebras
over~\(X\); this is the main point of
Section~\ref{sec:Cstar_over_X}.  Our definition is quite
natural, but more restrictive than the definitions in
\cites{Kirchberg:Michael, Bonkat:Thesis}.  The approach there
is to use the map \(\Open(X)\to\Open(\Prim A)\) induced by
\(\psi\colon \Prim(A)\to X\), where \(\Open(X)\) denotes the
complete lattice of open subsets of~\(X\).  If~\(X\) is a sober
space --~this is a very mild assumption that is also made under
a different name in \cites{Kirchberg:Michael,
  Bonkat:Thesis}~-- then we can recover it from the lattice
\(\Open(X)\), and a continuous map \(\Prim(A)\to X\) is
equivalent to a map \(\Open(X)\to\Open(\Prim A)\) that commutes
with arbitrary unions and finite intersections.

The definition of the Kasparov groups \(\KK_*(X;A,B)\) still
makes sense for \emph{any} map \(\Open(X)\to\Open(\Prim A)\)
(in the category of sets), that is, even the restrictions
imposed in \cites{Kirchberg:Michael, Bonkat:Thesis} can be
removed.  But such a map \(\Open(X)\to\Open(\Prim A)\)
corresponds to a continuous map \(\Prim(A)\to Y\) for another,
more complicated space~\(Y\) that contains~\(X\) as a subspace.
Hence the definitions in \cites{Kirchberg:Michael,
  Bonkat:Thesis} are, in fact, not more general.  But they
complicate computations because the discontinuities add further
input data which must be taken into account even for examples
where they vanish because the action is continuous.

Since the relevant point-set topology is widely unknown among
operator algebraists, we also recall some basic notions such as
sober spaces and Alexandrov spaces.  The latter are highly
non-Hausdorff spaces --~Alexandrov \(T_1\)\nb-spaces are
necessarily discrete~-- which are essentially the same as
preordered sets.  Any finite topological space is an Alexandrov
space, and their basic properties are crucial for this article.
To get acquainted with the setup, we simplify the description
of \(\Cst\)\nb-algebras over Alexandrov spaces and discuss some
small examples.  These rather elementary considerations
appeared previously in the theory of locales.

In Section~\ref{sec:KK_over_X}, we briefly recall the
definition and the basic properties of bivariant Kasparov
theory for \(\Cst\)\nb-algebras over a topological space.  We omit
most proofs because they are similar to the familiar arguments
for ordinary Kasparov theory and because the technical details
are already dealt with in~\cite{Bonkat:Thesis}.  We emphasise
the triangulated category structure on the Kasparov category
over~\(X\) because it plays an important role in connection
with the bootstrap class.

In Section~\ref{sec:bootstrap}, we define the bootstrap class
over a topological space~\(X\).  If~\(X\) is finite, we give
criteria for a \(\Cst\)\nb-algebra over~\(X\) to belong to the
bootstrap class.  These depend heavily on the relation between
Alexandrov spaces and preordered sets and therefore do not
extend directly to infinite spaces.

We define the \(X\)\nb-equivariant bootstrap class
\(\Bootstrap(X)\) as the localising subcategory of the Kasparov
category of \(\Cst\)\nb-algebras over~\(X\) that is generated by
the basic objects \((\C,x)\) for \(x\in X\), where we identify
\(x\in X\) with the corresponding constant map \(\Prim(\C)\to
X\).  Notice that this is exactly the list of all
\(\Cst\)\nb-algebras over~\(X\) with underlying
\(\Cst\)\nb-algebra~\(\C\).

We show that a \emph{nuclear} \(\Cst\)\nb-algebra \((A,\psi)\)
over~\(X\) belongs to the \(X\)\nb-equivariant bootstrap class
if and only if its ``fibres'' \(A(x)\) belong to the usual
bootstrap class for all \(x\in X\).  These fibres are certain
subquotients of~\(A\); if \(\psi\colon \Prim(A)\to X\) is a
homeomorphism, then they are exactly the \emph{simple}
subquotients of the \(\Cst\)\nb-algebra~\(A\).

The bootstrap class we define is the class of objects where we
expect a Universal Coefficient Theorem to hold.  If \(A\)
and~\(B\) belong to the bootstrap class, then an element of
\(\KK_*(X;A,B)\) is invertible if and only if it is fibrewise
invertible on \(\K\)\nb-theory, that is, the induced maps
\(\K_*\bigl(A(x)\bigr) \to \K_*\bigl(B(x)\bigr)\) are
invertible for all \(x\in X\).  This follows easily from our
definition of the bootstrap class.  The proof of our criterion
for a \(\Cst\)\nb-algebra over~\(X\) to belong to the bootstrap
class already provides a spectral sequence that computes
\(\KK_*(X;A,B)\) in terms of non-equivariant Kasparov groups.
Unfortunately, this spectral sequence is not useful for
classification purposes because it rarely degenerates to an
\emph{exact} sequence.

We call a \(\Cst\)\nb-algebra over~\(X\) \emph{tight} if the
map \(\Prim(A)\to X\) is a homeomorphism.  This implies that
its fibres are simple.  We show in
Section~\ref{sec:simple_fibres} that any separable nuclear
\(\Cst\)\nb-algebra over~\(X\) is \(\KK(X)\)\nb-equivalent to a
tight, separable, nuclear, purely infinite, stable
\(\Cst\)\nb-algebra over~\(X\).  The main issue is tightness.
By Kirchberg's classification result, this model is unique up
to \(X\)\nb-equivariant \Star{}isomorphism.  In this sense,
tight, separable, nuclear, purely infinite, stable
\(\Cst\)\nb-algebras over~\(X\) are classified up to
isomorphism by the isomorphism classes of objects in a certain
triangulated category: the subcategory of nuclear
\(\Cst\)\nb-algebras over~\(X\) in the Kasparov category.  The
difficulty is to replace this complete ``invariant'' by a more
tractable one that classifies objects of the --~possibly
smaller~-- bootstrap category \(\Bootstrap(X)\) by
\(\K\)\nb-theoretic data.

\smallskip

If~\(\Cat\) is a category, then we write \(A\inOb\Cat\) to
denote that~\(A\) is an object of~\(\Cat\) --~as opposed to a
morphism in~\(\Cat\).

\section{\texorpdfstring{$\Cst$}{C*}-algebras over a topological
  space}
\label{sec:Cstar_over_X}

We define the category \(\Cstarcat(X)\) of \(\Cst\)\nb-algebras
over a topological space~\(X\).  In the Hausdorff case, this
amounts to the familiar category of
\(\CONT_0(X)\)-\(\Cst\)\nb-algebras.  For non-Hausdorff spaces, our
notion is related to another one by Eberhard Kirchberg.  For
the Universal Coefficient Theorem, we must add some continuity
conditions to Kirchberg's definition of \(\Cstarcat(X)\).  We
explain in~\S\ref{sec:discontinuous_bundles} why these
conditions result in essentially no loss of generality.
Furthermore, we explain briefly why it is allowed to restrict
to the case where the underlying space~\(X\) is sober, and we
consider some examples, focusing on special properties of
finite spaces and Alexandrov spaces.

\subsection{The Hausdorff case}
\label{sec:Hausdorff_case}

Let~\(A\) be a \(\Cst\)\nb-algebra and let~\(X\) be a locally
compact Hausdorff space.  There are various equivalent
additional structures on~\(A\) that turn it into a
\(\Cst\)\nb-algebra over~\(X\) (see~\cite{Nilsen:Bundles} for the
proofs of most of the following assertions).  The most common
definition is the following one from~\cite{Kasparov:Novikov}:

\begin{definition}
  \label{def:Cstar_over_Hausdorff}
  A \emph{\(\CONT_0(X)\)-\(\Cst\)\nb-algebra} is a
  \(\Cst\)\nb-algebra~\(A\) together with an essential
  \Star{}homomorphism~\(\varphi\) from \(\CONT_0(X)\) to the
  centre of the multiplier algebra of~\(A\).  We abbreviate
  \(h\cdot a\defeq \varphi(h)\cdot a\) for \(h\in\CONT_0(X)\).

  A \Star{}homomorphism \(f\colon A\to B\) between two
  \(\CONT_0(X)\)-\(\Cst\)\nb-algebras is
  \emph{\(\CONT_0(X)\)-linear} if \(f(h\cdot a) = h\cdot f(a)\)
  for all \(h\in\CONT_0(X)\), \(a\in A\).

  Let \(\Cstarcat\bigl(\CONT_0(X)\bigr)\) be the category of
  \(\CONT_0(X)\)-\(\Cst\)\nb-algebras, whose morphisms are the
  \(\CONT_0(X)\)-linear \Star{}homomorphisms.
\end{definition}

A map~\(\varphi\) as above is equivalent to an \(A\)\nb-linear
essential \Star{}homomorphism
\[
\bar\varphi\colon \CONT_0(X,A)
\cong \CONT_0(X)\otimes_\maxi A\to A,
\qquad f\otimes a\mapsto \varphi(f)\cdot a,
\]
which exists by the universal property of the maximal tensor
product; the centrality of~\(\varphi\) ensures
that~\(\bar\varphi\) is a \Star{}homomorphism and well-defined.
Conversely, we get~\(\varphi\) back from~\(\bar\varphi\) by
restricting to elementary tensors; the assumed
\(A\)\nb-linearity of~\(\bar\varphi\) ensures that
\(\varphi(h)\cdot a \defeq \bar\varphi(h\otimes a)\) is a
multiplier of~\(A\).  The description via~\(\bar\varphi\) has
two advantages: it requires no multipliers, and the resulting
class in \(\KK_0(\CONT_0(X,A),A)\) plays a role in connection
with duality in bivariant Kasparov theory
(see~\cite{Emerson-Meyer:Dualities}).

Any \(\CONT_0(X)\)-\(\Cst\)\nb-algebra is isomorphic to the
\(\Cst\)\nb-algebra of \(\CONT_0(X)\)\nb-sections of an upper
semi-continuous \(\Cst\)\nb-algebra bundle over~\(X\)
(see~\cite{Nilsen:Bundles}).  Even more, this yields an
equivalence of categories between
\(\Cstarcat\bigl(\CONT_0(X)\bigr)\) and the category of upper
semi-continuous \(\Cst\)\nb-algebra bundles over~\(X\).

\begin{definition}
  \label{def:Prim}
  Let \(\Prim(A)\) denote the \emph{primitive ideal space}
  of~\(A\), equipped with the usual hull--kernel topology, also
  called Jacobson topology.
\end{definition}

The Dauns--Hofmann Theorem identifies the centre of the
multiplier algebra of~\(A\) with the \(\Cst\)\nb-algebra
\(\CONT_\bo\bigl(\Prim(A)\bigr)\) of bounded continuous
functions on the primitive ideal space of~\(A\).  Therefore,
the map~\(\varphi\) in
Definition~\ref{def:Cstar_over_Hausdorff} is of the form
\[
\psi^*\colon \CONT_0(X) \to \CONT_\bo(\Prim A),
\qquad f\mapsto f\circ\psi,
\]
for some continuous map \(\psi\colon \Prim(A)\to X\)
(see~\cite{Nilsen:Bundles}).  Thus \(\varphi\) and~\(\psi\) are
equivalent additional structures.  We use such maps~\(\psi\) to
generalise Definition~\ref{def:Cstar_over_Hausdorff} to the
non-Hausdorff case.

\subsection{The general definition}
\label{sec:Cstar_over_top_def}

Let~\(X\) be an arbitrary topological space.

\begin{definition}
  \label{def:Cstar_over_X}
  A \emph{\(\Cst\)\nb-algebra over~\(X\)} is a pair \((A,\psi)\)
  consisting of a \(\Cst\)\nb-algebra~\(A\) and a continuous map
  \(\psi\colon \Prim(A)\to X\).
\end{definition}

Our next task is to define morphisms between \(\Cst\)\nb-algebras
\(A\) and~\(B\) over the same space~\(X\).  This requires some
care because the primitive ideal space is \emph{not} functorial
for arbitrary \Star{}homomorphisms.

\begin{definition}
  \label{def:Open_X}
  For a topological space~\(X\), let \(\Open(X)\) be the set of
  open subsets of~\(X\), partially ordered by~\(\subseteq\).
\end{definition}

\begin{definition}
  \label{def:Ideals_A}
  For a \(\Cst\)\nb-algebra~\(A\), let \(\Ideals(A)\) be the set of
  all closed \Star{}ideals in~\(A\), partially ordered
  by~\(\subseteq\).
\end{definition}

The partially ordered sets \((\Open(X),\subseteq)\) and
\((\Ideals(A),\subseteq)\) are \emph{complete lattices}, that
is, any subset in them has both an infimum \(\bigwedge S\) and
a supremum \(\bigvee S\).  Namely, in \(\Open(X)\), the
supremum is~\(\bigcup S\), and the infimum is the
\emph{interior of~\(\bigcap S\)}; in \(\Ideals(A)\), the
infimum and supremum are
\[
\bigwedge_{I\in S} I = \bigcap_{I\in S} I,
\qquad
\bigvee_{I\in S} I = \cl{\sum_{I\in S} I}.
\]

We always identify \(\Open\bigl(\Prim(A)\bigr)\) and
\(\Ideals(A)\) using the isomorphism
\begin{equation}
  \label{eq:Prim_open_Ideal}
  \Open\bigl(\Prim(A)\bigr) \cong \Ideals(A),
  \qquad
  U \mapsto \bigcap_{\primid\in\Prim(A)\setminus U} \primid
\end{equation}
(see \cite{Dixmier:Cstar-algebres}*{\S3.2}).  This is a lattice
isomorphism and hence preserves infima and suprema.

Let \((A,\psi)\) be a \(\Cst\)\nb-algebra over~\(X\).  We get a map
\[
\psi^*\colon \Open(X)\to \Open(\Prim A) \cong \Ideals(A),
\qquad
U \mapsto \{\primid\in \Prim(A)\mid \psi(\primid)\in U\}
\cong A(U).
\]
We usually write \(A(U)\in\Ideals(A)\) for the ideal and
\(\psi^*(U)\) or \(\psi^{-1}(U)\) for the corresponding open
subset of \(\Prim(A)\).  If~\(X\) is a locally compact
Hausdorff space, then \(A(U)\defeq \CONT_0(U)\cdot A\) for all
\(U\in\Open(X)\).

\begin{example}
  \label{exa:Cstar_over_Prim}
  For any \(\Cst\)\nb-algebra~\(A\), the pair \((A,\ID_{\Prim A})\)
  is a \(\Cst\)\nb-algebra over \(\Prim(A)\); the ideals \(A(U)\)
  for \(U\in\Open(\Prim A)\) are given
  by~\eqref{eq:Prim_open_Ideal}.  \(\Cst\)\nb-algebras over
  topological spaces of this form play an important role
  in~\S\ref{sec:simple_fibres}, where we call them
  \emph{tight}.
\end{example}

\begin{lemma}
  \label{lem:psi_star_properties}
  The map~\(\psi^*\) is compatible with arbitrary suprema
  \textup{(}unions\textup{)} and finite infima
  \textup{(}intersections\textup{)}, so that
  \[
  A\biggl(\bigcup_{U\in S} U\biggr) = \cl{\sum_{U\in S} A(U)}, \qquad
  A\biggl(\bigcap_{U\in F} U\biggr) = \bigcap_{U\in F} A(U)
  \]
  for any subset \(S\subseteq\Open(X)\) and for any finite
  subset \(F\subseteq \Open(X)\).
\end{lemma}

\begin{proof}
  This is immediate from the definition.
\end{proof}

Taking for \(S\) and~\(F\) the empty set, this specialises to
\(A(\emptyset) = \{0\}\) and \(A(X)=A\).  Taking \(S=\{U,V\}\)
with \(U\subseteq V\), this specialises to the monotonicity
property
\[
U\subseteq V \qquad\Longrightarrow\qquad A(U)\subseteq A(V);
\]
We will implicitly use later that these properties follow from
compatibility with finite infima and suprema.

The following lemma clarifies when the map~\(\psi^*\) is
compatible with infinite infima.

\begin{lemma}
  \label{lem:open_map_infima}
  If the map \(\psi\colon \Prim(A)\to X\) is open or if~\(X\)
  is finite, then the map \(\psi^*\colon
  \Open(X)\to\Ideals(A)\) preserves infima --~that is, it maps
  the interior of \(\bigcap_{U\in S} U\) to the ideal
  \(\bigcap_{U\in S} A(U)\) for any subset
  \(S\subseteq\Open(X)\).  Conversely, if~\(\psi^*\) preserves
  infima and~\(X\) is a \(T_1\)\nb-space, that is, points
  in~\(X\) are closed, then~\(\psi\) is open.
\end{lemma}

Since preservation of infinite infima is automatic for
finite~\(X\), the converse assertion cannot hold for
general~\(X\).

\begin{proof}
  If~\(X\) is finite, then any subset of \(\Open(X)\) is
  finite, and there is nothing more to prove.  Suppose
  that~\(\psi\) is open.  Let~\(V\) be the interior of
  \(\bigcap_{U\in S} U\).  Let \(W\subseteq\Prim(A)\) be the
  open subset that corresponds to the ideal \(\bigcap_{U\in S}
  \psi^*(U)\).  We must show \(\psi^*(V)=W\).  Monotonicity
  yields \(\psi^*(V)\subseteq W\).  Since~\(\psi\) is open,
  \(\psi(W)\) is an open subset of~\(X\).  By construction,
  \(\psi(W)\subseteq U\) for all \(U\in S\) and hence
  \(\psi(W)\subseteq V\).  Thus \(\psi^*(V)\supseteq
  \psi^*\bigl(\psi(W)\bigr) \supseteq W\supseteq \psi^*(V)\),
  so that \(\psi^*(V)=W\).

  Now suppose, conversely, that~\(\psi^*\) preserves infima and
  that points in~\(X\) are closed.  Assume that~\(\psi\) is not
  open.  Then there is an open subset~\(W\) in \(\Prim(A)\) for
  which \(\psi(W)\) is not open in~\(X\).  Let \(S\defeq \{
  X\setminus\{x\} \mid x\in X\setminus \psi(W)\}\subseteq
  \Open(X)\); this is where we need points to be closed.  We
  have \(\bigcap_{U\in S} U = \psi(W)\) and \(\bigcap_{U\in S}
  \psi^*(U) = \psi^{-1}\bigl(\psi(W)\bigr)\).  Since
  \(\psi(W)\) is not open, the infimum~\(V\) of~\(S\) in
  \(\Open(X)\) is strictly smaller than \(\psi(W)\).  Hence
  \(\psi^*(V)\) cannot contain~\(W\).  But~\(W\) is an open
  subset of \(\psi^{-1}\bigl(\psi(W)\bigr)\) and hence
  contained in the infimum of \(\psi^*(S)\) in \(\Open(\Prim
  A)\).  Therefore, \(\psi^*\) does not preserve infima,
  contrary to our assumption.  Hence~\(\psi\) must be open.
\end{proof}

For a locally compact Hausdorff space~\(X\), the map
\(\Prim(A)\to X\) is open if and only if~\(A\) corresponds to a
\emph{continuous} \(\Cst\)\nb-algebra bundle over~\(X\)
(see \cite{Nilsen:Bundles}*{Theorem 2.3}).

\begin{definition}
  \label{def:morphism_over_X}
  Let \(A\) and~\(B\) be \(\Cst\)\nb-algebras over a topological
  space~\(X\).  A \Star{}homomorphism \(f\colon A\to B\) is
  \emph{\(X\)\nb-equivariant} if \(f\bigl(A(U)\bigr)\subseteq
  B(U)\) for all \(U\in\Open(X)\).
\end{definition}

For locally compact Hausdorff spaces, this is equivalent to
\(\CONT_0(X)\)-linearity by the following variant of
\cite{Bonkat:Thesis}*{Propositon 5.4.7}:

\begin{proposition}
  \label{pro:equivariance_linearity}
  Let \(A\) and~\(B\) be \(\Cst\)\nb-algebras over a locally
  compact Hausdorff space~\(X\), and let \(f\colon A\to B\) be
  a \Star{}homomorphism.  The following assertions are
  equivalent:
  \begin{enumerate}[label=\textup{(\arabic{*})}]
  \item \(f\) is \(\CONT_0(X)\)-linear;

  \item \(f\) is \(X\)\nb-equivariant, that is,
    \(f\bigl(A(U)\bigr)\subseteq B(U)\) for all
    \(U\in\Open(X)\);

  \item \(f\) descends to the fibres, that is,
    \(f\bigl(A(X\setminus\{x\})\bigr) \subseteq
    B(X\setminus\{x\})\) for all \(x\in X\).
  \end{enumerate}
\end{proposition}

To understand the last condition, recall that the fibres of the
\(\Cst\)\nb-algebra bundle associated to~\(A\) are \(A_x\defeq A
\mathbin/ A(X\setminus\{x\})\).  Condition~(3) means that~\(f\)
descends to maps \(f_x\colon A_x\to B_x\) for all \(x\in X\).

\begin{proof}
  It is clear that
  (1)\(\Longrightarrow\)(2)\(\Longrightarrow\)(3).  The
  equivalence (3)\(\iff\)(1) is the assertion of
  \cite{Bonkat:Thesis}*{Propositon 5.4.7}.  To check that~(3)
  implies~(1), take \(h\in\CONT_0(X)\) and \(a\in A\).  We get
  \(f(h\cdot a)= h\cdot f(a)\) provided both sides have the
  same values at all \(x\in X\) because the map \(A\to
  \prod_{x\in X} A_x\) is injective.  Now~(3) implies
  \(f(h\cdot a)_x = h(x)\cdot f(a)_x\) because
  \(\bigl(h-h(x)\bigr)\cdot a\in A(X\setminus\{x\})\).
\end{proof}

\begin{definition}
  \label{def:Cstar_over_X_cat}
  Let \(\Cstarcat(X)\) be the category whose objects are the
  \(\Cst\)\nb-algebras over~\(X\) and whose morphisms are the
  \(X\)\nb-equivariant \Star{}homomorphisms.  We write
  \(\Hom_X(A,B)\) for this set of morphisms.
\end{definition}

Proposition~\ref{pro:equivariance_linearity} yields an
isomorphism of categories \(\Cstarcat\bigl(\CONT_0(X)\bigr)
\cong \Cstarcat(X)\).  In this sense, our theory for general
spaces extends the more familiar theory of
\(\CONT_0(X)\)-\(\Cst\)\nb-algebras.

\subsection{Locally closed subsets and subquotients}
\label{sec:subquotients}

\begin{definition}
  \label{def:loclo}
  A subset~\(C\) of a topological space~\(X\) is called
  \emph{locally closed} if it is the intersection of an open
  and a closed subset or, equivalently, of the form
  \(C=U\setminus V\) with \(U,V\in\Open(X)\); we can also
  assume \(V\subseteq U\) here.  We let \(\Loclo(X)\) be the
  set of locally closed subsets of~\(X\).
\end{definition}

A subset is locally closed if and only if it is relatively open
in its closure.  Being locally closed is inherited by finite
intersections, but \emph{not} by unions or complements.

\begin{definition}
  \label{def:loclo_subquotients}
  Let~\(X\) be a topological space and let \((A,\psi)\) be a
  \(\Cst\)\nb-algebra over~\(X\).  Write \(C\in\Loclo(X)\) as
  \(C=U\setminus V\) for open subsets \(U,V\subseteq X\) with
  \(V\subseteq U\).  We define
  \[
  A(C) \defeq A(U) \mathbin/ A(V).
  \]
\end{definition}

\begin{lemma}
  \label{lem:loclo_subquotients}
  The subquotient \(A(C)\) does not depend on \(U\) and~\(V\)
  above.
\end{lemma}

\begin{proof}
  Let \(U_1,V_1,U_2,V_2\in\Open(X)\) satisfy \(V_1\subseteq
  U_1\), \(V_2\subseteq U_2\), and \(U_1\setminus
  V_1=U_2\setminus V_2\).  Then \(V_1\cup U_2 = U_1\cup U_2 =
  U_1\cup V_2\) and \(V_1\cap U_2=V_1\cap V_2=U_1\cap V_2\).
  Since \(U\mapsto A(U)\) preserves unions, this implies
  \[
  A(U_2)+ A(V_1)=A(U_1)+A(V_2).
  \]
  We divide this equation by \(A(V_1\cup V_2)=A(V_1)+A(V_2)\).
  This yields
  \[
  \frac{A(U_2)+ A(V_1)}{A(V_1\cup V_2)}
  \cong \frac{A(U_2)}{A(U_2) \cap A(V_1\cup V_2)}
  = \frac{A(U_2)}{A\bigl(U_2 \cap (V_1\cup V_2)\bigr)}
  = \frac{A(U_2)}{A(V_2)}
  \]
  on the left hand side and, similarly, \(A(U_1)\mathbin/
  A(V_1)\) on the right hand side.  Hence \(A(U_1)\mathbin/
  A(V_1) \cong A(U_2)\mathbin/ A(V_2)\) as desired.
\end{proof}

Now assume that \(X=\Prim(A)\) and \(\psi=\ID_{\Prim(A)}\).
Lemma~\ref{lem:loclo_subquotients} associates a subquotient
\(A(C)\) of~\(A\) to each locally closed subset of
\(\Prim(A)\).  Equation~\eqref{eq:Prim_open_Ideal} shows that
any subquotient of~\(A\) arises in this fashion; here
subquotient means: a quotient of one ideal in~\(A\) by another
ideal in~\(A\).  Open subsets of~\(X\) correspond to ideals,
closed subsets to quotients of~\(A\).  For any
\(C\in\Loclo(\Prim A)\), there is a canonical homeomorphism
\(\Prim\bigl(A(C)\bigr) \cong C\).  This is well-known if~\(C\)
is open or closed, and the general case reduces to these
special cases.

\begin{example}
  \label{exa:simple_subquotients}
  If \(\Prim(A)\) is a finite topological \(T_0\)\nb-space,
  then any singleton~\(\{\primid\}\) in \(\Prim(A)\) is locally
  closed (this holds more generally for the Alexandrov
  \(T_0\)\nb-spaces introduced in~\S\ref{sec:top_poset} and
  follows from the description of closed subsets in terms of
  the specialisation preorder).

  Since \(\Prim A(C)\cong C\), the subquotients
  \(A_\primid\defeq A(\{\primid\})\) for \(\primid\in\Prim(A)\)
  are precisely the \emph{simple subquotients} of~\(A\).
\end{example}

\begin{example}
  \label{exa:bad_interval}
  Consider the interval \([0,1]\) with the topology where the
  non-empty closed subsets are the closed intervals \([a,1]\)
  for all \(a\in[0,1]\).  A non-empty subset is locally closed
  if and only if it is either of the form \([a,1]\) or
  \([a,b)\) for \(a,b\in[0,1]\) with \(a<b\).  In this space,
  singletons are not locally closed.  Hence a
  \(\Cst\)\nb-algebra with this primitive ideal space has no
  simple subquotients.
\end{example}

\subsection{Functoriality and tensor products}
\label{sec:Cstar_over_top_functor}

\begin{definition}
  \label{def:Cstarcat_X_functor}
  Let \(X\) and~\(Y\) be topological spaces.  A continuous map
  \(f\colon X\to Y\) induces a functor
  \[
  f_*\colon \Cstarcat(X) \to \Cstarcat(Y),
  \qquad
  (A,\psi)\mapsto (A,f\circ\psi).
  \]
  Thus \(X\mapsto \Cstarcat(X)\) is a functor from the category
  of topological spaces to the category of categories (up to
  the usual issues with sets and classes).
\end{definition}

Since \((f\circ\psi)^{-1}=\psi^{-1}\circ f^{-1}\), we have
\[
(f_*A)(C) = A\bigl(f^{-1}(C)\bigr)
\qquad
\text{for all \(C\in\Loclo(Y)\).}
\]

If \(f\colon X\to Y\) is the embedding of a subset with the
subspace topology, we also write
\[
i_X^Y\defeq f_*\colon \Cstarcat(X)\to\Cstarcat(Y)
\]
and call this the \emph{extension} functor from~\(X\) to~\(Y\).
We have \((i_X^YA)(C) = A(C\cap X)\) for all \(C\in\Loclo(Y)\).

\begin{definition}
  \label{def:restrict_bundle}
  Let~\(X\) be a topological space and let~\(Y\) be a
  \emph{locally closed} subset of~\(X\), equipped with the
  subspace topology.  Let \((A,\psi)\) be a \(\Cst\)\nb-algebra
  over~\(X\).  Its \emph{restriction to~\(Y\)} is a
  \(\Cst\)\nb-algebra~\(A|_Y\) over~\(Y\), consisting of the
  \(\Cst\)\nb-algebra \(A(Y)\) defined as in
  Definition~\ref{def:loclo_subquotients}, equipped with the
  canonical map
  \[
  \Prim A(Y) \congto \psi^{-1}(Y) \xrightarrow{\psi} Y.
  \]
  Thus \(A|_Y(C) = A(C)\) for
  \(C\in\Loclo(Y)\subseteq\Loclo(X)\).
\end{definition}

It is clear that the restriction to~\(Y\) provides a functor
\[
r_X^Y\colon \Cstarcat(X)\to\Cstarcat(Y)
\]
that satisfies \(r_Y^Z \circ r_X^Y =r_X^Z\) if \(Z\subseteq
Y\subseteq X\) and \(r_X^X=\ID\).

If \(Y\) and~\(X\) are Hausdorff and locally compact, then a
continuous map \(f\colon Y\to X\) also induces a pull-back
functor
\begin{multline*}
  f^*\colon \Cstarcat(X) \cong \Cstarcat\bigl(\CONT_0(X)\bigr)
  \to \Cstarcat\bigl(\CONT_0(Y)\bigr) \cong \Cstarcat(Y),\\
  A\mapsto \CONT_0(Y)\otimes_{\CONT_0(X)} A.
\end{multline*}
For the constant map \(Y\to\pt\), this functor
\(\Cstarcat\to\Cstarcat(Y)\) maps a \(\Cst\)\nb-algebra~\(A\) to
\(f^*(A)\defeq \CONT_0(Y,A)\) with the obvious
\(\CONT_0(Y)\)-\(\Cst\)\nb-algebra structure.  This functor has no
analogue for a non-Hausdorff space~\(Y\).  Therefore, a
continuous map \(f\colon Y\to X\) need not induce a functor
\(f^*\colon \Cstarcat(X)\to\Cstarcat(Y)\).  For embeddings of
locally closed subsets, the functor~\(r_X^Y\) plays the role
of~\(f^*\).

\begin{lemma}
  \label{lem:adjointness}
  Let~\(X\) be a topological space and let \(Y\subseteq X\).

  \begin{enumerate}[label=\textup{(\alph*)}]
  \item If~\(Y\) is \emph{open}, then there are natural
    isomorphisms
    \[
    \Hom_X(i_Y^X(A),B) \cong \Hom_Y\bigl(A,r_X^Y(B)\bigr)
    \]
    if \(A\) and~\(B\) are \(\Cst\)\nb-algebras over \(Y\)
    and~\(X\), respectively.

    In other words, \(i_Y^X\) is left adjoint to~\(r_X^Y\).

  \item If~\(Y\) is \emph{closed}, then there are natural
    isomorphisms
    \[
    \Hom_Y(r_X^Y(A),B) \cong \Hom_X\bigl(A,i_Y^X(B)\bigr)
    \]
    if \(A\) and~\(B\) are \(\Cst\)\nb-algebras over \(X\)
    and~\(Y\), respectively.

    In other words, \(i_Y^X\) is right adjoint to~\(r_X^Y\).

  \item For any locally closed subset \(Y\subseteq X\), we have
    \(r_X^Y\circ i_Y^X(A)=A\) for all \(\Cst\)\nb-algebras~\(A\)
    over~\(Y\).

  \end{enumerate}
\end{lemma}

\begin{proof}
  We first prove~(a).  We have \(i_Y^X(A)(U) = A(U\cap Y)\) for
  all \(U\in\Open(X)\), and this is an ideal in \(A(U)\).  A
  morphism \(\varphi\colon i_Y^X(A)\to B\) is equivalent to a
  \Star{}homomorphism \(\varphi\colon A(Y)\to B(X)\) that maps
  \(A(U\cap Y)\to B(U)\) for all \(U\in\Open(X)\).  This holds
  for all \(U\in\Open(X)\) once it holds for \(U\in\Open(Y)
  \subseteq \Open(X)\).  Hence~\(\varphi\) is equivalent to a
  \Star{}homomorphism \(\varphi'\colon A(Y)\to B(Y)\) that maps
  \(A(U)\to B(U)\) for all \(U\in\Open(Y)\).  The latter is
  nothing but a morphism \(A\to r_X^Y(B)\).  This proves~(a).

  Now we turn to~(b).  Again, we have \(i_Y^X(B)(U) = B(U\cap
  Y)\) for all \(U\in\Open(X)\), but now this is a quotient of
  \(B(U)\).  A morphism \(\varphi\colon A\to i_Y^X(B)\) is
  equivalent to a \Star{}homomorphism \(\varphi\colon A(X)\to
  B(Y)\) that maps \(A(U)\to B(U\cap Y)\) for all
  \(U\in\Open(X)\).  Hence \(A(X\setminus Y)\) is mapped to
  \(B(\emptyset)=0\), so that~\(\varphi\) descends to a
  map~\(\varphi'\) from \(A \mathbin/ A(X\setminus Y)\cong
  A(Y)\) to \(B(Y)\) that maps \(A(U\cap Y)\) to~\(B(U)\) for
  all \(U\in\Open(X)\).  The latter is equivalent to a morphism
  \(r_X^Y(A)\to B\) as desired.  This finishes the proof
  of~(b).

  Assertion~(c) is trivial.
\end{proof}

\begin{example}
  \label{exa:extend_from_point}
  For each \(x\in X\), we get a map \(i_x=i_x^X\colon \pt\cong
  \{x\} \subseteq X\) from the one-point space to~\(X\).  The
  resulting functor \(\Cstarcat \to \Cstarcat(X)\) maps a
  \(\Cst\)\nb-algebra~\(A\) to the \(\Cst\)\nb-algebra \(i_x(A)=(A,x)\)
  over~\(X\), where~\(x\) also denotes the constant map
  \[
  x\colon \Prim(A)\to X,
  \qquad \primid\mapsto x\quad
  \text{for all \(\primid\in\Prim(A)\).}
  \]
  If \(C\in\Loclo(X)\), then
  \[
  i_x(A)(C) =
  \begin{cases}
    A & \text{if \(x\in C\);}\\
    0 & \text{otherwise.}
  \end{cases}
  \]
  The functor~\(i_x\) plays an important role if~\(X\) is
  finite.  The generators of the bootstrap class are of the
  form \(i_x(\C)\).  Each \(\Cst\)\nb-algebra over~\(X\) carries a
  canonical filtration whose subquotients are of the form
  \(i_x(A)\).
\end{example}

\begin{lemma}
  \label{lem:Hom_ix}
  Let~\(X\) be a topological space and let \(x\in X\).  Then
  \[
  \Hom_X\bigl(A,i_x^X(B)\bigr)
  \cong \Hom\bigl(A\bigl(\cl{\{x\}}\bigr),B\bigr)
  \]
  for all \(A\inOb\Cstarcat(X)\), \(B\inOb\Cstarcat\), and
  \[
  \Hom_X(i_x^X(A),B)
  \cong \Hom\Bigl(A,\bigcap_{U\in \mathcal{U}_x} B(U)\Bigr).
  \]
  for all \(A\inOb\Cstarcat\), \(B\inOb\Cstarcat(X)\),
  where~\(\mathcal{U}_x\) denotes the open neighbourhood filter
  of~\(x\) in~\(X\).  If~\(x\) has a minimal open
  neighbourhood~\(U_x\), then this becomes
  \[
  \Hom_X(i_x^X(A),B)
  \cong \Hom\bigl(A,B(U_x)\bigr).
  \]
\end{lemma}

Recall that \(A\inOb\Cat\) means that~\(A\) is an object
of~\(\Cat\).

\begin{proof}
  Let \(C\defeq \cl{\{x\}}\).  Then any non-empty open subset
  \(V\subseteq C\) contains~\(x\), so that \(i_x^C(B)(V)=B\).
  This implies \(\Hom_C(A,i_x^C(B)) \cong
  \Hom\bigl(A(C),B\bigr)\).  Combining this with
  \(i_x^X=i_C^X\circ i_x^C\) and the adjointness relation in
  Lemma~\ref{lem:adjointness}.(b) yields
  \[
  \Hom_X\bigl(A,i_x^X(B)\bigr)
  \cong \Hom_C\bigl(r_X^C(A), i_x^C(B)\bigr)
  \cong \Hom(A(C),B).
  \]

  An \(X\)\nb-equivariant \Star{}homomorphism \(i_x^XA\to B\)
  restricts to a family of compatible maps \(A = (i_x^XA)(U)\to
  B(U)\) for all \(U\in\mathcal{U}_x\), so that we get a
  \Star{}homomorphism from~\(A\) to
  \(\bigcap_{U\in\mathcal{U}_x} B(U)\).  Conversely, any such
  \Star{}homomorphism \(A\to \bigcap_{U\in\mathcal{U}_x} B(U)\)
  provides an \(X\)\nb-equivariant \Star{}homomorphism
  \(i_x^XA\to B\).  This yields the second assertion.
\end{proof}

Let \(A\) and~\(B\) be \(\Cst\)\nb-algebras and let \(A\minotimes
B\) be their minimal (or spatial) \(\Cst\)\nb-tensor product.  Then
there is a canonical continuous map
\[
\Prim(A) \times \Prim(B)\to \Prim(A\minotimes B).
\]
Therefore, if \(A\) and~\(B\) are \(\Cst\)\nb-algebras over \(X\)
and~\(Y\), respectively, then \(A\minotimes B\) is a
\(\Cst\)\nb-algebra over \(X\times Y\).  This defines a bifunctor
\[
{\minotimes}\colon \Cstarcat(X)\times\Cstarcat(Y)\to
\Cstarcat(X\times Y).
\]

In particular, if \(Y=\pt\) is the one-point space, then we get
endofunctors \(\blank\otimes B\) on \(\Cstarcat(X)\) for
\(B\inOb\Cstarcat\) because \(X\times\pt\cong X\).

If~\(X\) is a Hausdorff space, then the diagonal in \(X\times
X\) is closed and we get an internal tensor product
functor~\(\otimes_X\) in \(\Cstarcat(X)\) by restricting the
external tensor product in \(\Cstarcat(X\times X)\) to the
diagonal.  This operation has no analogue for general~\(X\).

\subsection{Restriction to sober spaces}
\label{sec:sober}

A space is sober if and only if it can be recovered from its
lattice of open subsets.  Any topological space can be
completed to a sober space with the same lattice of open
subsets.  Therefore, it usually suffices to study
\(\Cst\)\nb-algebras over sober topological spaces.

\begin{definition}
  \label{def:sober}
  A topological space is \emph{sober} if each irreducible
  closed subset of~\(X\) is the closure \(\cl{\{x\}}\) of
  exactly one singleton of~\(X\).  Here an \emph{irreducible}
  closed subset of~\(X\) is a non-empty closed subset of~\(X\)
  which is not the union of two proper closed subsets of
  itself.
\end{definition}

If~\(X\) is not sober, let~\(\hat{X}\) be the set of all
irreducible closed subsets of~\(X\).  There is a canonical map
\(\iota\colon X\to\hat{X}\) which sends a point \(x\in X\) to
its closure.  If \(S\subseteq X\) is closed, let
\(\hat{S}\subseteq\hat{X}\) be the set of all \(A\in\hat{X}\)
with \(A\subseteq S\).  The map \(S\mapsto\hat{S}\) commutes
with finite unions and arbitrary intersections; in particular,
it maps~\(X\) itself to all of~\(\hat{X}\) and~\(\emptyset\) to
\(\hat{\emptyset}=\emptyset\).  Hence the subsets
of~\(\hat{X}\) of the form~\(\hat{S}\) for closed subsets
\(S\subseteq X\) form the closed subsets of a topology
on~\(\hat{X}\).

The map~\(\iota\) induces a bijection between the families of
closed subsets of \(X\) and~\(\hat{X}\).  Hence~\(\iota\) is
continuous, and it induces a bijection \(\iota^*\colon
\Open(\hat{X})\to\Open(X)\).  It also follows that~\(\hat{X}\)
is a sober space because \(X\) and~\(\hat{X}\) have the same
irreducible closed subsets.

Since the morphisms in \(\Cstarcat(X)\) only use \(\Open(X)\),
the functor
\[
\iota_*\colon \Cstarcat(X) \to \Cstarcat(\hat{X})
\]
is fully faithful.  Therefore, we do not lose much if we assume
our topological spaces to be sober.

The following example shows a pathology that can occur if the
separation axiom~\(T_0\) fails:

\begin{example}
  \label{exa:chaotic_topology}
  Let~\(X\) carry the chaotic topology
  \(\Open(X)=\{\emptyset,X\}\).  Then \(\hat{X}=\star\) is the
  space with one point.  By definition, an action of~\(X\) on a
  \(\Cst\)\nb-algebra~\(A\) is a map \(\Prim(A)\to X\).  But for a
  \Star{}homomorphism \(A\to B\) between two \(\Cst\)\nb-algebras
  over~\(X\), the \(X\)\nb-equivariance condition imposes no
  restriction.  Hence all maps \(\Prim(A)\to X\) yield
  isomorphic objects of \(\Cstarcat(X)\).
\end{example}

\begin{lemma}
  \label{lem:sober_open_to_space}
  If~\(X\) is a sober topological space, then there is a
  bijective correspondence between continuous maps
  \(\Prim(A)\to X\) and maps \(\Open(X)\to\Ideals(A)\) that
  commute with arbitrary suprema and finite infima; it sends a
  continuous map \(\psi\colon \Prim(A)\to X\) to the map
  \[
  \psi^*\colon
  \Open(X)\to\Open\bigl(\Prim(A)\bigr)=\Ideals(A).
  \]
\end{lemma}

\begin{proof}
  We have already seen that a continuous map
  \(\psi\colon\Prim(A)\to X\) generates a map~\(\psi^*\) with
  the required properties for any space~\(X\).

  Conversely, let \(\psi^*\colon \Open(X)\to \Ideals(A)\) be a
  map that preserves arbitrary unions and finite intersections.
  Given \(\primid\in\Prim(A)\), let~\(U_\primid\) be the union
  of all \(U\in\Open(X)\) with \(\primid\notin\psi^*(U)\).
  Then \(\primid\notin\psi^*(U_\primid)\) because~\(\psi^*\)
  preserves unions, and~\(U_\primid\) is the maximal open
  subset with this property.  Thus \(A_\primid\defeq X\setminus
  U_\primid\) is the minimal closed subset with
  \(\primid\notin\psi^*(X\setminus A_\primid)\).  This subset
  is non-empty because \(\psi^*(X) = \Prim(A)\)
  contains~\(\primid\), and irreducible because~\(\psi^*\)
  preserves finite intersections.

  Since~\(X\) is sober, there is a unique \(\psi(\primid)\in
  X\) with \(A_\primid = \cl{\{\psi(\primid)\}}\).  This
  defines a map \(\psi\colon \Prim(A)\to X\).  If \(U\subseteq
  X\) is open, then \(\psi(\primid)\notin U\) if and only if
  \(A_\primid\cap U=\emptyset\), if and only if
  \(\primid\notin\psi^*(U)\).  Hence \(\psi^*(U) =
  \psi^{-1}(U)\).  This shows that~\(\psi\) is continuous and
  generates~\(\psi^*\).  Thus the map \(\psi\to\psi^*\) is
  surjective.

  Since sober spaces are~\(T_0\), two different continuous maps
  \(\psi_1,\psi_2\colon \Prim(A)\to X\) generate different maps
  \(\psi_1^*,\psi_2^*\colon \Open(X)\to \Ideals(A)\).  Hence
  the map \(\psi\to\psi^*\) is also injective.
\end{proof}

\subsection{Some very easy examples}
\label{sec:Cstar_over_X_easy_examples}

Here we describe the categories of \(\Cst\)\nb-algebras over the
three sober topological spaces with at most two points.

\begin{example}
  \label{exa:Cstar_over_point}
  If~\(X\) is a single point, then \(\Cstarcat(X)\) is
  isomorphic to the category of \(\Cst\)\nb-algebras (without any
  extra structure).
\end{example}

Up to homeomorphism, there are two sober topological spaces
with two points.  The first one is the discrete space.

\begin{example}
  \label{exa:Cstar_over_two-point}
  The category of \(\Cst\)\nb-algebras over the discrete two-point
  space is equivalent to the product category
  \(\Cstarcat\times\Cstarcat\) of pairs of \(\Cst\)\nb-algebras.
\end{example}

More generally, if \(X=X_1\sqcup X_2\) is a disjoint union of
two subspaces, then
\begin{equation}
  \label{eq:disjoint_product}
  \Cstarcat(X) \simeq \Cstarcat(X_1)\times \Cstarcat(X_2).
\end{equation}
Thus it usually suffices to study connected spaces.

\begin{example}
  \label{exa:extensions}
  Another sober topological space with two points is
  \(X=\{1,2\}\) with
  \[
  \Open(X) = \bigl\{\emptyset,\{1\},\{1,2\}\bigr\}.
  \]
  A \(\Cst\)\nb-algebra over this space comes with a single
  distinguished ideal \(A(1)\triangleleft A\), which is
  arbitrary.  Thus we get the category of pairs \((I,A)\)
  where~\(I\) is an ideal in~\(A\).  We may associate to this
  data the \(\Cst\)\nb-algebra extension \(I\into A\prto A/I\).  In
  fact, the morphisms in \(\Hom_X(A,B)\) are the morphisms of
  extensions
  \[\xymatrix{
    A(1)\ \ar[d] \ar@{>->}[r]&
    A \ar[d] \ar@{->>}[r]&
    A\mathbin/ A(1) \ar[d]\\
    B(1)\ \ar@{>->}[r]&
    B \ar@{->>}[r]&
    B\mathbin/ B(1).
  }
  \]
  Thus \(\Cstarcat(X)\) is equivalent to the category of
  \emph{\(\Cst\)\nb-algebra extensions}.  This example is also
  studied in~\cite{Bonkat:Thesis}.
\end{example}

\subsection{Topologies and partial orders}
\label{sec:top_poset}

Certain non-Hausdorff spaces are closely related to partially
ordered sets.  In particular, there is a bijection between
sober topologies and partial orders on a finite set.  Here we
recall the relevant constructions.

\begin{definition}
  \label{def:specialisation_preorder}
  Let~\(X\) be a topological space.  The \emph{specialisation
    preorder}~\(\preceq\) on~\(X\) is defined by \(x\preceq y\)
  if the closure of~\(\{x\}\) is contained in the closure
  of~\(\{y\}\) or, equivalently, if~\(y\) is contained in all
  open subsets of~\(X\) that contain~\(x\).  Two points \(x\)
  and~\(y\) are called \emph{topologically indistinguishable}
  if \(x\preceq y\) and \(y\preceq x\), that is, the closures of
  \(\{x\}\) and~\(\{y\}\) are equal.
\end{definition}

The separation axiom~\(T_0\) means that topologically
indistinguishable points are equal.  Since this is automatic
for sober spaces, \(\preceq\) is a partial order on~\(X\) in
all cases we need.  As usual, we write \(x\prec y\) if
\(x\preceq y\) and \(x\neq y\), and \(x\succeq y\) and \(x\succ
y\) are equivalent to \(y\preceq x\) and \(y\prec x\),
respectively.

The separation axiom~\(T_1\) requires points to be closed.
This is equivalent to the partial order~\(\preceq\) being
trivial, that is, \(x\preceq y\) if and only if \(x=y\).  Thus
our partial order is only meaningful for highly non-separated
spaces.

The following notion goes back to an article by Paul Alexandrov
from 1937 (\cite{Alexandrov:Diskrete_Raeume}); see
also~\cite{Arenas:Alexandrov} for a more recent reference, or
the english
\href{http://en.wikipedia.org/wiki/Alexandrov_topology}{Wikipedia}
entry on the Alexandrov topology.

\begin{definition}
  \label{def:Alexandrov_open}
  Let \((X,\le)\) be a preordered set.  A subset \(S\subseteq
  X\) is called \emph{Alexandrov-open} if \(S\ni x\le y\)
  implies \(y\in S\).  The Alexandrov-open subsets form a
  topology on~\(X\) called the \emph{Alexandrov topology}.
\end{definition}

A subset of~\(X\) is closed in the Alexandrov topology if and
only if \(S\ni x\) and \(x\ge y\) imply \(S\ni y\).  It is
locally closed if and only if it is \emph{convex}, that is,
\(x\le y\le z\) and \(x,z\in S\) imply \(y\in S\).  In
particular, singletons are locally closed (compare
Example~\ref{exa:simple_subquotients}).

The specialisation preorder for the Alexandrov topology is the
given preorder.  Moreover, a map \((X,\le)\to (Y,\le)\) is
continuous for the Alexandrov topology if and only if it is
monotone.  Thus we have identified the category of preordered
sets with monotone maps with a full subcategory of the category
of topological spaces.

It also follows that if a topological space carries an
Alexandrov topology for some preorder, then this preorder must
be the specialisation preorder.  In this case, we call the
space an \emph{Alexandrov space} or a \emph{finitely generated
  space}.  The following lemma provides some equivalent
descriptions of Alexandrov spaces; the last two explain in what
sense these spaces are finitely generated.

\begin{lemma}
  \label{lem:Alexandrov_spaces}
  Let~\(X\) be a topological space.  The following are
  equivalent:
  \begin{itemize}
  \item \(X\) is an Alexandrov space;
  \item an arbitrary intersection of open subsets of~\(X\) is
    open;

  \item an arbitrary union of closed subsets of~\(X\) is
    closed;

  \item every point of~\(X\) has a smallest neighbourhood;

  \item a point~\(x\) lies in the closure of a subset~\(S\)
    of~\(X\) if and only if \(x\in\cl{\{y\}}\) for some \(y\in
    S\);

  \item \(X\) is the inductive limit of the inductive system of
    its finite subspaces.

  \end{itemize}
\end{lemma}

\begin{corollary}
  \label{cor:finite_Alexandrov}
  Any finite topological space is an Alexandrov space.  Thus
  the construction of Alexandrov topologies and specialisation
  preorders provides a bijection between preorders and
  topologies on a finite set.
\end{corollary}

\begin{definition}
  \label{def:minimal_open}
  Let~\(X\) be an Alexandrov space.  We denote the minimal open
  neighbourhood of \(x\in X\) by \(U_x\in\Open(X)\).
\end{definition}

We have
\[
U_x\subseteq U_y
\iff x\in U_y
\iff y\in \cl{\{x\}}
\iff \cl{\{y\}}\subseteq \cl{\{x\}}
\iff y\preceq x.
\]

If~\(X\) is an Alexandrov space, then we can simplify the data
for a \(\Cst\)\nb-algebra over~\(X\) as follows:

\begin{lemma}
  \label{lem:simplify_Cstar_over_finite}
  A \(\Cst\)\nb-algebra over a sober Alexandrov space~\(X\) is
  determined uniquely by a \(\Cst\)\nb-algebra~\(A\) together with
  ideals \(A(U_x)\triangleleft A\) for all \(x\in X\), subject
  to the two conditions \(\cl{\sum_{x\in X} A(U_x)} = A\) and
  \begin{equation}
    \label{eq:compatibility_AUX}
    A(U_x)\cap A(U_y) = \cl{\sum_{z\in U_x\cap U_y} A(U_z)}
    \qquad \text{for all \(x,y\in X\).}
  \end{equation}
\end{lemma}

\begin{proof}
  A map \(\Open(X)\to\Ideals(A)\) that preserves suprema and
  maps~\(U_x\) to \(A(U_x)\) for all \(x\in X\) must map \(U =
  \bigvee_{x\in U} U_x\) to \(\bigvee_{x\in U} A(U_x) =
  \cl{\sum_{x\in U} A(U_x)}\).  The map so defined preserves
  suprema by construction.  The two hypotheses of the lemma
  ensure \(A(X)=A\) and \(A(U_x\cap U_y) = A(U_x)\cap A(U_y)\)
  for all \(x,y\in X\).  Hence they are necessary for
  preservation of finite infima.

  Since the lattice \(\Ideals(A)\cong \Open(\Prim A)\) is
  distributive, \eqref{eq:compatibility_AUX} implies
  \begin{multline*}
    A(U) \wedge A(V)
    = \bigvee_{x\in U} A(U_x) \wedge \bigvee_{y\in V} A(V_y)
    = \bigvee_{(x,y)\in U\times V} A(U_x) \wedge A(V_y)
    \\= \bigvee_{(x,y)\in U\times V} A(U_x\cap V_y)
    = A(U\cap V);
  \end{multline*}
  the last step uses that \(U\mapsto A(U)\) commutes with
  suprema.  We clearly have \(A(\emptyset)=\{0\}\) as well, so
  that \(U\mapsto A(U)\) preserves arbitrary finite
  intersections.  Therefore, our map \(\Open(X)\to\Ideals(A)\)
  satisfies the conditions in
  Lemma~\ref{lem:sober_open_to_space} and hence comes from a
  continuous map \(\Prim A\to X\).
\end{proof}

Of course, a \Star{}homomorphism \(A\to B\) between two
\(\Cst\)\nb-algebras over~\(X\) is \(X\)\nb-equivariant if and
only if it maps \(A(U_x)\to B(U_x)\) for all \(x\in X\).

Equation~\eqref{eq:compatibility_AUX} implies \(A(U_x)\subseteq
A(U_y)\) if \(U_x\subseteq U_y\), that is, if \(x\succeq y\).
Thus the map \(x\mapsto A(U_x)\) is order-reversing.  It
sometimes happens that \(U_x\cap U_y=U_z\) for some \(x,y,z\in
X\).  In this case, we may drop the ideal \(A(U_z)\) from the
description of a \(\Cst\)\nb-algebra over~\(X\) and replace the
condition~\eqref{eq:compatibility_AUX} for \(x,y\) by
\(A(U_w)\subseteq A(U_x)\cap A(U_y)\) for all \(w\in U_x\cap
U_y\).

\subsection{Some more examples}
\label{sec:more_examples}

A useful way to represent finite partially ordered sets and
hence finite sober topological spaces is via finite
\emph{directed acyclic graphs}.

To a partial order~\(\preceq\) on~\(X\), we associate the
finite directed acyclic graph with vertex set~\(X\) and with an
arrow \(x\leftarrow y\) if and only if \(x\prec y\) and there
is no \(z\in X\) with \(x\prec z\prec y\).  We can recover the
partial order from this graph by letting \(x\preceq y\) if and
only if the graph contains a directed path \(x\leftarrow
x_1\leftarrow \dotsb \leftarrow x_n\leftarrow y\).

We have reversed arrows here because an arrow \(x\to y\) means
that \(A(U_x)\subseteq A(U_y)\).  Furthermore, \(x\in U_y\) if
and only if there is a directed path from~\(x\) to~\(y\).  Thus
we can read the meaning of the
relations~\eqref{eq:compatibility_AUX} from the graph.

\begin{example}
  \label{exa:totally_ordered}
  Let \((X,\ge)\) be a set with a total order, such as
  \(\{1,\dotsc,n\}\) with the order~\(\ge\).  The corresponding
  graph is
  \[
  \xymatrix{1\ar[r]&2\ar[r]&3\ar[r]&\dotsb\ar[r]&n.}
  \]
  For totally ordered~\(X\), \eqref{eq:compatibility_AUX} is
  equivalent to monotonicity of the map \(x\mapsto A(U_x)\)
  with respect to the opposite order~\(\le\) on~\(X\).  As a
  consequence, a \(\Cst\)\nb-algebra over~\(X\) is nothing but
  a \(\Cst\)\nb-algebra~\(A\) together with a monotone map
  \((X,\le)\to\Ideals(A)\), \(x\mapsto A(U_x)\), such that
  \(\bigvee_{x\in X} A(U_x)=A\).  For \(X=
  \bigl(\{1,\dotsc,n\},\ge\bigr)\), the latter condition just
  means \(A(U_n)=A\), so that we can drop this ideal.  Thus we
  get \(\Cst\)\nb-algebras with an increasing chain of \(n-1\)
  ideals \(I_1\triangleleft I_2\triangleleft\dotsb\triangleleft
  I_{n-1} \triangleleft A\).  This situation is studied in
  detail in~\cite{Meyer-Nest:Filtrated_K}.
\end{example}

\begin{figure}
  \begin{gather*}
    \xymatrix{\bullet\ar[r]&\bullet\ar[r]&\bullet}\qquad
    \xymatrix{\bullet\ar[r]\ar[dr]&\bullet\\&\bullet}\qquad
    \xymatrix{\bullet\ar[r]&\bullet\\\bullet\ar[ur]}\\
    \xymatrix{\bullet\ar[r]&\bullet\ar[r]&\bullet\ar[r]&\bullet} \qquad
    \xymatrix{\bullet\ar[r]&\bullet\ar[r]\ar[dr]&\bullet\\&&\bullet} \qquad
    \xymatrix{\bullet\ar[r]\ar[dr]&\bullet\ar[r]&\bullet\\&\bullet}\\
    \xymatrix{\bullet\ar[r]\ar[d]&\bullet\ar[d]\\\bullet\ar[r]&\bullet} \qquad
    \xymatrix{\bullet\ar[r]\ar[dr]&\bullet\\\bullet\ar[ur]\ar[r]&\bullet} \qquad
    \xymatrix{\bullet\ar[r]\ar[dr]&\bullet\\\bullet\ar[r]&\bullet} \qquad
    \xymatrix{\bullet\ar[r]&\bullet\ar[r]&\bullet\\\bullet\ar[ur]} \\
    \xymatrix{&\bullet\\\bullet\ar[r]\ar[dr]\ar[ur]&\bullet\\&\bullet} \qquad
    \xymatrix{\bullet\ar[dr]&\\\bullet\ar[r]&\bullet\\\bullet\ar[ur]&}
  \end{gather*}
  \caption{Connected directed acyclic graphs with three or four vertices}
  \label{fig:graphs}
\end{figure}
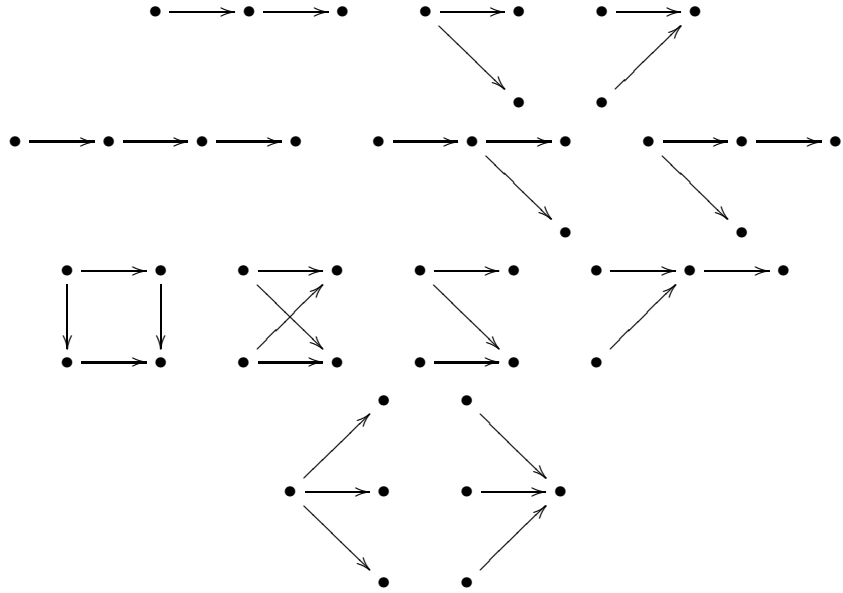

Using that any finite topological space is an Alexandrov space,
we can easily list all homeomorphism classes of finite
topological spaces with, say, three or four elements.  We only
consider sober spaces here, and we assume connectedness to
further reduce the number of cases.  Under these assumptions,
Figure~\ref{fig:graphs} contains a complete list.
The first and fourth case are already contained in
Example~\ref{exa:totally_ordered}.
Lemma~\ref{lem:simplify_Cstar_over_finite} describes
\(\Cst\)\nb-algebras over the spaces in Figure~\ref{fig:graphs}
as \(\Cst\)\nb-algebras equipped with three or four ideals
\(A(U_x)\) for \(x\in X\), subject to some conditions, which
often make some of the ideals redundant.

\begin{example}
  \label{exa:second_case}
  The second graph in Figure~\ref{fig:graphs} describes
  \(\Cst\)\nb-algebras with three ideals \(A(U_j)\), \(j=1,2,3\),
  subject to the conditions \(A(U_2)\cap A(U_3) = A(U_1)\) and
  \(A(U_2)+A(U_3)=A\).  This is equivalent to prescribing only
  two ideals \(A(U_2)\) and \(A(U_3)\) subject to the single
  condition \(A(U_2)+A(U_3)=A\).
\end{example}

\begin{example}
  \label{exa:third_case}
  Similarly, the third graph in Figure~\ref{fig:graphs}
  describes \(\Cst\)\nb-algebras with two distinguished ideals
  \(A(U_1)\) and \(A(U_2)\) subject to the condition
  \(A(U_1)\cap A(U_2) = \{0\}\); here \(U_3=X\) implies
  \(A(U_3)=A\).
\end{example}

\begin{example}
  \label{exa:ninth_case}
  The ninth case above is more complicated.  We label our
  points by \(1,2,3,4\) such that \(1\rightarrow 3\leftarrow
  2\rightarrow 4\).  Here we have a \(\Cst\)\nb-algebra~\(A\) with
  four ideals \(I_j\defeq A(U_j)\) for \(j=1,2,3,4\), subject
  to the conditions
  \[
  I_1\subseteq I_3,\qquad
  I_1\cap I_4=\{0\},\qquad
  I_2= I_3\cap I_4,\qquad
  I_3+I_4 = A.
  \]
  Thus the ideal~\(I_2\) is redundant, and we are left with
  three ideals \(I_1,I_3,I_4\) subject to the conditions
  \(I_1\subseteq I_3\), \(I_1\cap I_4=\{0\}\), and \(I_3+I_4 =
  A\).
\end{example}

\subsection{How to treat discontinuous bundles}
\label{sec:discontinuous_bundles}

The construction of \(X\)\nb-equivariant Kasparov theory in
\cites{Kirchberg:Michael,Bonkat:Thesis} works for \emph{any}
map \(\psi^*\colon \Open(X)\to\Ideals(A)\), we do not need the
conditions in Lemma~\ref{lem:sober_open_to_space}.  Here we
show how to reduce this more general situation to the case
considered above: discontinuous actions of \(\Open(X)\) as in
\cites{Kirchberg:Michael,Bonkat:Thesis} are equivalent to
continuous actions of another space~\(Y\) that contains~\(X\)
as a subspace.  The category \(\Cstarcat(Y)\)
contains~\(\Cstarcat(X)\) as a full subcategory, and a similar
statement holds for the associated Kasparov categories.  As a
result, allowing general maps~\(\psi^*\) merely amounts to
replacing the space~\(X\) by the larger space~\(Y\).  For
\(\Cst\)\nb-algebras that really live over the subspace~\(X\), the
extension to~\(Y\) significantly complicates the computation of
the Kasparov groups.  This is why we always require~\(\psi^*\)
to satisfy the conditions in
Lemma~\ref{lem:sober_open_to_space}, which ensure that it comes
from a continuous map \(\Prim(A)\to X\).

\begin{example}
  \label{exa:discrete_discontinuous_action}
  Let \(X=\{1,2\}\) with the discrete topology.  A monotone map
  \(\psi^*\colon \Open(X)\to A\) with
  \(\psi^*(\emptyset)=\{0\}\) and \(\psi^*(X)=A\) as considered
  in \cites{Bonkat:Thesis,Kirchberg:Michael} is equivalent to
  specifying two \emph{arbitrary} ideals \(A(1)\) and \(A(2)\).
  This automatically generates the ideals \(A(1)\cap A(2)\) and
  \(A(1)\cup A(2)\).  We can encode these four ideals in an
  action of a topological space~\(Y\) with four points
  \(\{1\cap 2,1,2,3\}\) and open subsets
  \[
  \emptyset,\quad
  \{1\cap 2\},\quad
  \{1\cap 2,1\},\quad
  \{1\cap 2,2\},\quad
  \{1\cap 2,1,2\},\quad
  \{1\cap 2,1,2,3\}.
  \]
  The corresponding graph is the seventh one in
  Figure~\ref{fig:graphs}.  The map~\(\psi^*\) maps these open
  subsets to the ideals
  \[
  \{0\},\quad
  A(1)\cap A(2),\quad
  A(1),\quad
  A(2),\quad
  A(1)\cup A(2),\quad
  A,
  \]
  respectively.  This defines a complete lattice morphism
  \(\Open(Y)\to\Ideals(A)\), and any complete lattice morphism
  is of this form for two ideals \(A(1)\) and \(A(2)\).  Thus
  an action of \(\Open(\{1,2\})\) in the generalised sense
  considered in \cites{Bonkat:Thesis,Kirchberg:Michael} is
  equivalent to an action of~\(Y\) in our sense.

  Any \(X\)\nb-equivariant \Star{}homomorphism \(A\to B\)
  between two such discontinuous \(\Cst\)\nb-algebras over~\(X\)
  will also preserve the ideals \(A(1)\cap A(2)\) and
  \(A(1)\cup A(2)\).  Hence it is \(Y\)\nb-equivariant as well.
  Therefore, the above construction provides an equivalence of
  categories between \(\Cstarcat(Y)\) and the category of
  \(\Cst\)\nb-algebras with an action of \(\Open(X)\) in the sense
  of \cites{Bonkat:Thesis,Kirchberg:Michael}.

  Whereas the computation of
  \[
  \KK_*(X;A,B) \cong \KK_*\bigl(A(1),B(1)\bigr)\times
  \KK_*\bigl(A(2),B(2)\bigr)
  \]
  for two \(\Cst\)\nb-algebras \(A\) and~\(B\) over~\(X\) is
  trivial, the corresponding problem for \(\Cst\)\nb-algebras
  over~\(Y\) is an interesting problem: this is one of the
  small examples where filtrated \(\K\)\nb-theory does not yet
  suffice for classification.
\end{example}

This simple example generalises as follows.  Let \(f\colon
\Open(X)\to \Ideals(A) \cong \Open(\Prim A)\) be an arbitrary
map.  Let \(Y\defeq 2^{\Open(X)}\) be the power set of
\(\Open(X)\), partially ordered by inclusion.  We describe the
topology on~\(Y\) below.  We embed the original space~\(X\)
into~\(Y\) by mapping \(x\in X\) to its open neighbourhood
filter:
\[
\mathcal{U}\colon X\to Y,
\qquad
x\mapsto \{U\in\Open(X)\mid x\in U\}.
\]
We define a map
\[
\psi\colon \Prim(A)\to Y,
\qquad \primid\mapsto \{U\in\Open(X)\mid \primid\in f(U)\}.
\]
For \(y\in Y\), let \(Y_{\supseteq y}\defeq \{x\in Y\mid
x\supseteq y\}\).  For a singleton \(\{U\}\) with
\(U\in\Open(X)\), we easily compute
\[
\psi^{-1}(Y_{\supseteq \{U\}}) = f(U)
\in \Ideals(A) \cong \Open(\Prim A).
\]
Moreover, \(Y_{\supseteq y\cup z} = Y_{\supseteq y} \cap
Y_{\supseteq z}\), so that we get
\[
\psi^{-1}(Y_{\supseteq \{U_1,\dotsc,U_n\}}) = f(U_1)
\cap\dotsb\cap f(U_n).
\]
A similar argument shows that
\[
\mathcal{U}^{-1}(Y_{\supseteq \{U_1,\dotsc,U_n\}})
= \mathcal{U}^{-1}(Y_{\supseteq U_1}) \cap\dotsb\cap
\mathcal{U}^{-1}(Y_{\supseteq U_n})
= U_1\cap \dotsb\cap U_n.
\]

We equip~\(Y\) with the topology that has the sets
\(Y_{\supseteq F}\) for \emph{finite} subsets~\(F\)
of~\(\Open(X)\) as a basis.  It is clear from the above
computations that this makes the maps \(\psi\)
and~\(\mathcal{U}\) continuous; even more, the subspace
topology on the range of~\(\mathcal{U}\) is the given topology
on~\(X\).

As a consequence, any map \(f\colon \Open(X)\to\Ideals(A)\)
turns~\(A\) into a \(\Cst\)\nb-algebra over the space
\(Y\supseteq X\).  Conversely, given a \(\Cst\)\nb-algebra
over~\(Y\), we define \(f\colon \Open(X)\to\Ideals(A)\) by
\(f(U) \defeq \psi^{-1}(Y_{\supseteq\{U\}})\).  This
construction is inverse to the one above.  Furthermore, a
\Star{}homomorphism \(A\to B\) that maps \(f_A(U)\) to
\(f_B(U)\) for all \(U\in\Open(X)\) also maps \(\psi^*_A(U)\)
to \(\psi^*_B(U)\) for all \(U\in\Open(Y)\).  We can sum this
up as follows:

\begin{theorem}
  \label{the:discontinuous_actions}
  The category of \(\Cst\)\nb-algebras equipped with a map
  \(f\colon \Open(X)\to\Ideals(A)\) is isomorphic to the
  category of \(\Cst\)\nb-algebras over~\(Y\).
\end{theorem}

If~\(f\) has some additional properties like monotonicity, or
is a lattice morphism, then this limits the range of the
map~\(\psi\) above and thus allows us to replace~\(Y\) by a
smaller subset.  In \cite{Kirchberg:Michael}*{Definition 1.3}
and \cite{Bonkat:Thesis}*{Definition 5.6.2}, an action of a
space~\(X\) on a \(\Cst\)\nb-algebra is defined to be a map
\(f\colon \Open(X)\to\Ideals(A)\) that is monotone and
satisfies \(f(\emptyset)=\{0\}\) and \(f(X)=A\).  These
assumptions are equivalent to
\[
U\in\psi(\primid),\quad U\subseteq V
\qquad\Longrightarrow\qquad V\in\psi(\primid)
\]
and \(\emptyset\notin\psi(\primid)\) and \(X\in\psi(\primid)\)
for all \(\primid\in\Prim(A)\).  Hence the category of
\(\Cst\)\nb-algebras with an action of~\(X\) in the sense
of~\cites{Bonkat:Thesis, Kirchberg:Michael} is equivalent to
the category of \(\Cst\)\nb-algebras over the space
\[
Y' \defeq \{y\subseteq \Open(X) \mid
y\ni U\subseteq V \Longrightarrow
V\in y,\ \emptyset\notin y,\ X\in y\},
\]
equipped with the subspace topology from~\(Y\).

\section{Bivariant  \texorpdfstring{$K$}{K}-theory for
  \texorpdfstring{$\Cst$}{C*}-algebras over topological spaces}
\label{sec:KK_over_X}

Let~\(X\) be a topological space.  Eberhard
Kirchberg~\cite{Kirchberg:Michael} and Alexander
Bonkat~\cite{Bonkat:Thesis} define Kasparov groups
\(\KK_*(X;A,B)\) for separable \(\Cst\)\nb-algebras \(A\)
and~\(B\) over~\(X\).  More precisely, instead of a continuous
map \(\Prim(A)\to X\) they use a separable
\(\Cst\)\nb-algebra~\(A\) with a monotone map \(\psi^*\colon
\Open(X)\to \Ideals(A)\) with \(A(\emptyset)=\{0\}\) and
\(A(X)=A\).  This is more general because any continuous map
\(\Prim(A)\to X\) generates such a map \(\psi^*\colon
\Open(X)\to \Ideals(A)\).  Hence their definitions apply to
\(\Cst\)\nb-algebras over~\(X\) in our sense.  We have
explained in~\S\ref{sec:discontinuous_bundles} why the setting
in~\cites{Kirchberg:Michael, Bonkat:Thesis} is, despite
appearences, not more general than our setting.

If~\(X\) is Hausdorff and locally compact, \(\KK_*(X;A,B)\)
agrees with Gennadi Kasparov's theory
\(\mathcal{R}\KK_*(X;A,B)\) defined in~\cite{Kasparov:Novikov}.
In this section, we recall the definition and some basic
properties of the functor \(\KK_*(X;A,B)\) and the resulting
category \(\KKcat(X)\), and we equip the latter with a
triangulated category structure.

\subsection{The definition}
\label{sec:def_KK}

We assume from now on that the topology on~\(X\) has a
countable basis, and we restrict attention to separable
\(\Cst\)\nb-algebras.

\begin{definition}
  \label{def:Cstar_sep}
  A \(\Cst\)\nb-algebra \((A,\psi)\) over~\(X\) is called
  \emph{separable} if~\(A\) is a separable \(\Cst\)\nb-algebra.
  Let \(\Cstarsep(X)\subseteq \Cstarcat(X)\) be the full
  subcategory of separable \(\Cst\)\nb-algebras over~\(X\).
\end{definition}

To describe the cycles for \(\KK_*(X;A,B)\) recall that the
usual Kasparov cyles for \(\KK_*(A,B)\) are of the form
\((\varphi,\Hils_B,F,\gamma)\) in the even case (for \(KK_0\))
and \((\varphi,\Hils_B,F)\) in the odd case (for \(KK_1\)),
where
\begin{itemize}
\item \(\Hils_B\) is a right Hilbert \(B\)\nb-module;

\item \(\varphi\colon A\to \Bound(\Hils_B)\) is a
  \Star{}representation;

\item \(F\in \Bound(\Hils_B)\);

\item \(\varphi(a)(F^2-1)\), \(\varphi(a)(F-F^*)\), and
  \([\varphi(a),F]\) are compact for all \(a\in A\);

\item in the even case, \(\gamma\) is a \(\Z/2\)-grading
  on~\(\Hils_B\) --~that is, \(\gamma^2=1\) and
  \(\gamma=\gamma^*\)~-- that commutes with \(\varphi(A)\) and
  anti-commutes with~\(F\).
\end{itemize}

The following definition of \(X\)\nb-equivariant bivariant
\(\K\)\nb-theory is equivalent to the ones in
\cites{Kirchberg:Michael, Bonkat:Thesis}, see
\cite{Kirchberg:Michael}*{Definition 4.1}, and
\cite{Bonkat:Thesis}*{Definition 5.6.11 and Satz 5.6.12}.

\begin{definition}
  \label{def:KK_cycle_over_X}
  Let \(A\) and~\(B\) be \(\Cst\)\nb-algebras over~\(X\) (or,
  more generally, \(\Cst\)\nb-algebras with a map
  \(\Open(X)\to\Ideals(A)\)).  A Kasparov cycle
  \((\varphi,\Hils_B,F,\gamma)\) or \((\varphi,\Hils_B,F)\) for
  \(\KK_*(A,B)\) is called \emph{\(X\)\nb-equivariant} if
  \[
  \varphi\bigl(A(U)\bigr)\cdot \Hils_B
  \subseteq \Hils_B\cdot B(U)\qquad
  \text{for all \(U\in\Open(X)\).}
  \]
  Let \(\KK_*(X;A,B)\) be the group of homotopy classes of such
  \(X\)\nb-equivariant Kasparov cycles for \(\KK_*(A,B)\); a
  homotopy is an \(X\)\nb-equivariant Kasparov cycle for
  \(\KK_*(A,C([0,1])\otimes B)\), where we view
  \(C([0,1])\otimes B\) as a \(\Cst\)\nb-algebra over~\(X\) in
  the usual way (compare~\S\ref{sec:Cstar_over_top_functor}).
\end{definition}

The subset \(\Hils_B\cdot B(U)\subseteq \Hils_B\) is a closed
linear subspace by the Cohen--Hewitt Factorisation Theorem.

If~\(X\) is Hausdorff, then the extra condition in
Definition~\ref{def:KK_cycle_over_X} is equivalent to
\(\CONT_0(X)\)-linearity of~\(\varphi\) (compare
Proposition~\ref{pro:equivariance_linearity}).  Thus the above
definition of \(\KK_*(X;A,B)\) agrees with the more familiar
definition of \(\mathcal{R}\KK_*(X;A,B)\)
in~\cite{Kasparov:Novikov}.

If \(X=\pt\) is the one-point space, the \(X\)\nb-equivariance
condition is empty and we get the plain Kasparov theory
\(\KK_*(\pt;A,B)=\KK_*(A,B)\).

The same arguments as usual show that \(\KK_*(X;A,B)\) remains
unchanged if we strengthen the conditions for Kasparov cycles
by requiring \(F=F^*\) and \(F^2=1\).

\subsection{Basic properties}
\label{sec:basic_properties}

The Kasparov theory defined above has all the properties that
we can expect from a bivariant \(\K\)\nb-theory.

\begin{enumerate}
\item The groups \(\KK_*(X;A,B)\) define a bifunctor from
  \(\Cstarsep(X)\) to the category of \(\Z/2\)-graded Abelian
  groups, contravariant in the first and covariant in the
  second variable.

\item There is a natural, associative \emph{Kasparov
    composition product}
  \[
  \KK_i(X;A,B)\times \KK_j (X; B,C)\to \KK_{i+j}(X;A,C)
  \]
  if \(A,B,C\) are \(\Cst\)\nb-algebras over~\(X\).

  Furthermore, there is a natural \emph{exterior product}
  \[
  \KK_i(X;A,B) \times \KK_j(Y;C,D) \to
  \KK_{i+j}(X\times Y;A\minotimes C,B\minotimes D)
  \]
  for two spaces \(X\) and~\(Y\) and \(\Cst\)\nb-algebras \(A\),
  \(B\) over~\(X\) and \(C\), \(D\) over~\(Y\).

  The existence and properties of the Kasparov composition
  product and the exterior product are verified in a more
  general context in \cite{Bonkat:Thesis}*{\S3.2}.
\end{enumerate}

\begin{definition}
  \label{def:KKcat}
  Let \(\KKcat(X)\) be the category whose objects are the
  separable \(\Cst\)\nb-algebras over~\(X\) and whose morphism sets
  are \(\KK_0(X;A,B)\).
\end{definition}

\begin{enumerate}[resume]
\item\label{zero-object} The zero \(\Cst\)\nb-algebra acts as a
\emph{zero object} in \(\KKcat(X)\), that is,
  \[
  \KK_*(X;\{0\},A) = 0 = \KK_*(X;A,\{0\})
  \qquad\text{for all \(A\inOb\KKcat(X)\).}
  \]

\item\label{coproduct} The \(\CONT_0\)\nb-direct sum of a sequence of
  \(\Cst\)\nb-algebras behaves like a \emph{coproduct}, that is,
  \[
  \KK_*\biggl(X;\bigoplus_{n\in\N} A_n,B\biggr)
  \cong \prod_{n\in\N} \KK_*(X;A_n,B)
  \]
  if \(A_n,B\inOb\KKcat(X)\) for all \(n\in\N\).

\item\label{product} The direct sum \(A\oplus B\) of two separable
  \(\Cst\)\nb-algebras \(A\) and~\(B\) over~\(X\) is a \emph{direct
    product} in \(\KKcat(X)\), that is,
  \[
  \KK_*(X;D,A\oplus B) \cong \KK_*(X;D,A)\oplus \KK_*(X;D,B)
  \]
  for all \(D\inOb\KKcat(X)\) (see \cite{Bonkat:Thesis}*{Lemma
    3.1.9}).
\end{enumerate}

Properties (\ref{zero-object})--(\ref{product}) are summarised
as follows:

\begin{proposition}
  \label{pro:KKcat_additive}
  The category \(\KKcat(X)\) is additive and has countable
  coproducts.
\end{proposition}

\begin{enumerate}[resume]
\item The exterior product is compatible with the Kasparov
  product, \(\CONT_0\)\nb-direct sums, and addition, that is,
  it defines a countably additive bifunctor
  \[
  \minotimes\colon \KKcat(X)\otimes \KKcat(Y) \to \KKcat(X\times Y).
  \]
  This operation is evidently associative.

\item\label{tensored} In particular, \(\KKcat(X)\) is tensored over
  \(\KKcat(\pt)\cong\KKcat\), that is, \(\otimes\) provides an
  associative bifunctor
  \[
  \minotimes\colon \KKcat(X)\otimes \KKcat \to \KKcat(X).
  \]

\item The bifunctor \((A,B)\mapsto \KK_*(X;A,B)\) satisfies
  Bott periodicity, homotopy invariance, and \(\Cst\)\nb-stability
  in each variable.  This follows from the corresponding
  properties of \(\KKcat\) using the tensor structure
  in~(\ref{tensored}).

  For instance, the Bott periodicity isomorphism
  \(\CONT_0(\R^2)\cong\C\) in \(\KKcat\) yields \(A\otimes
  \CONT_0(\R^2) \cong A\otimes \C \cong A\) in \(\KKcat(X)\)
  for all \(A\inOb\KKcat(X)\).

\item\label{functoriality} The functor \(f_*\colon
  \Cstarcat(X)\to \Cstarcat(Y)\) for a continuous map \(f\colon
  X\to Y\) descends to a functor
  \[
  f_*\colon \KKcat(X)\to\KKcat(Y).
  \]
  In particular, this covers the extension functors \(i_X^Y\)
  for a subspace \(X\subseteq Y\).

\item\label{restriction_functor} The restriction functor
  \(r_X^Y\) for \(Y\in\Loclo(X)\) also descends to a functor
  \[
  r_X^Y\colon \KKcat(X)\to\KKcat(Y).
  \]
\end{enumerate}

\begin{definition}[see \cite{Bonkat:Thesis}*{Definition 5.6.6}]
  \label{def:semi-split_ext}
  A diagram \(I\to E\to Q\) in \(\Cstarcat(X)\) is an
  \emph{extension} if, for all \(U\in\Open(X)\), the diagrams
  \(I(U)\to E(U)\to Q(U)\) are extensions of \(\Cst\)\nb-algebras.
  We write \(I\into E\prto Q\) to denote extensions.

  An extension is called \emph{split} if it splits by an
  \(X\)\nb-equivariant \Star{}homomorphism.

  An extension is called \emph{semi-split} if there is a
  completely positive, contractive section \(Q\to E\) that is
  \emph{\(X\)\nb-equivariant}, that is, it restricts to
  sections \(Q(U)\to E(U)\) for all \(U\in\Open(X)\).
\end{definition}

If \(I\into E\prto Q\) is an extension of \(\Cst\)\nb-algebras
over~\(X\), then we get \(\Cst\)\nb-algebra extensions \(I(Y)\into
E(Y)\prto Q(Y)\) for all locally closed subsets \(Y\subseteq
X\).  If the original extension is semi-split, so are the
extensions \(I(Y)\into E(Y)\prto Q(Y)\) for \(Y\in\Loclo(X)\).
Even more, the functor \(r_X^Y\colon
\Cstarcat(X)\to\Cstarcat(Y)\) maps extensions in
\(\Cstarcat(X)\) to extensions in \(\Cstarcat(Y)\), and
similarly for split and semi-split extensions.

\begin{theorem}
  \label{the:excision}
  Let \(I\into E\prto Q\) be a semi-split extension in
  \(\Cstarsep(X)\) and let~\(B\) be a separable \(\Cst\)\nb-algebra
  over~\(X\).  There are six-term exact sequences
  \begin{gather*}
    \xymatrix{
      \KK_0(X;Q,B)\ar[r]&
      \KK_0(X;E,B)\ar[r]&
      \KK_0(X;I,B)\ar[d]^\partial\\
      \KK_1(X;I,B)\ar[u]^\partial&
      \KK_1(X;E,B)\ar[l]&
      \KK_1(X;Q,B)\ar[l]
    }\\\shortintertext{and}
    \xymatrix{
      \KK_0(X;B,I)\ar[r]&
      \KK_0(X;B,E)\ar[r]&
      \KK_0(X;B,Q)\ar[d]^\partial\\
      \KK_1(X;B,Q)\ar[u]^\partial&
      \KK_1(X;B,E)\ar[l]&
      \KK_1(X;B,I),\ar[l]
    }
  \end{gather*}
  where the horizontal maps in both exact sequences are induced
  by the given maps \(I\to E\to Q\), and the vertical maps are,
  up to signs, Kasparov products with the class of our
  semi-split extension in \(\KK_1(Q,I)\).

  Furthermore, extensions with a completely positive section
  are semi-split.
\end{theorem}

\begin{proof}
  The long exact sequences for semi-split extensions follow
  from \cite{Bonkat:Thesis}*{Satz 3.3.10} or from
  \cite{Bonkat:Thesis}*{Korollar 5.6.13}.

  The last sentence plays a technical role in the proof of
  Proposition~\ref{pro:nuclear_local}.  We have to replace an
  \(X\)\nb-equivariant completely positive section \(s\colon
  Q\to E\) by another section that is an \(X\)\nb-equivarant
  completely positive contraction.  Without
  \(X\)\nb-equivariance, this is done in
  \cite{Cuntz-Skandalis:Cones}*{Remark 2.5}.  We claim that the
  constructions during the proof yield \(X\)\nb-equivariant
  maps if we start with \(X\)\nb-equivariant maps.

  Let \(Q^+\) and~\(E^+\) be obtained by adjoining units to
  \(Q\) and~\(E\).  Let \((u_n)_{n\in\N}\) be an approximate
  unit in~\(Q\) and let \(v_n \defeq \sup(1,s(u_n))\)
  in~\(E^+\).  Since \(v_n\ge1\), \(v_n\) is invertible.  The maps
  \[
  s_n\colon Q^+\to E^+,
  \qquad q\mapsto v_n^{-\nicefrac{1}{2}}
  s(u_n^{\nicefrac{1}{2}} q u_n^{\nicefrac{1}{2}})
  v_n^{-\nicefrac{1}{2}}
  \]
  considered in~\cite{Cuntz-Skandalis:Cones} are unital and
  completely positive and hence contractive.  They are
  \(X\)\nb-equivariant if~\(s\) is \(X\)\nb-equivariant.
  Since~\(v_n\) lifts \(1\in Q^+\), the maps~\((s_n)\) lift the
  maps \(q\mapsto u_n^{\nicefrac{1}{2}} q
  u_n^{\nicefrac{1}{2}}\), which converge pointwise to the
  identity map.

  It remains to show that the space of maps \(Q^+\to Q^+\) that
  lift to an \(X\)\nb-equivariant unital completely positive
  map \(Q^+\to E^+\) is closed in the topology of pointwise
  norm convergence.  Without \(X\)\nb-equivariance, this is
  \cite{Arveson:Extensions_Cstar}*{Theorem 6}.  Its proof is
  based on the following construction.  Let
  \(\varphi,\psi\colon Q^+\rightrightarrows E^+\) be two unital
  completely positive maps and let~\((e_m)_{m\in\N}\) be a
  quasi-central approximate unit for~\(I\) in~\(Q^+\).  Then
  Arveson uses the unital completely positive maps
  \[
  q\mapsto
  e_m^{\nicefrac{1}{2}} \varphi(q) e_m^{\nicefrac{1}{2}}
  + (1-e_m)^{\nicefrac{1}{2}} \psi(q) (1-e_m)^{\nicefrac{1}{2}}.
  \]
  Clearly, this map is \(X\)\nb-equivariant, if \(\varphi\)
  and~\(\psi\) are \(X\)\nb-equivariant.  Hence the argument
  in~\cite{Arveson:Extensions_Cstar} produces a Cauchy sequence
  of \(X\)\nb-equivariant unital completely positive maps
  \(\hat{s}_n\colon Q^+\to E^+\) lifting~\(s_n\).  Its limit is
  an \(X\)\nb-equivariant unital completely positive section
  \(Q^+\to E^+\).
\end{proof}

\begin{theorem}
  \label{the:KKX_universal}
  The canonical functor \(\Cstarsep(X)\to\KKcat(X)\) is the
  universal split-exact \(\Cst\)\nb-stable
  \textup{(}homotopy\textup{)} functor.
\end{theorem}

\begin{proof}
  This follows from \cite{Bonkat:Thesis}*{Satz 3.5.10}, compare
  also \cite{Bonkat:Thesis}*{Korollar 5.6.13}.  The homotopy
  invariance assumption is redundant because, by a deep theorem
  of Nigel Higson, a split-exact, \(\Cst\)\nb-stable functor is
  automatically homotopy invariant.  This holds for
  \(\Cstarsep\) itself and is inherited by \(\Cstarsep(X)\)
  because of the tensor product operation
  \(\Cstarsep(X)\times\Cstarsep\to\Cstarsep(X)\).
\end{proof}

\subsection{Triangulated category structure}
\label{sec:tri_KKcatX}

We are going to turn \(\KKcat(X)\) into a triangulated category
as in~\cite{Meyer-Nest:BC}.  We have already remarked that
\(\KKcat(X)\) is additive.  The \emph{suspension functor} is
\(\Sigma(A)\defeq \CONT_0(\R,A) = \CONT_0(\R)\minotimes A\).
This functor is an automorphism (up to natural isomorphisms) by
Bott periodicity.

The \emph{mapping cone triangle}
\begin{equation}
  \label{eq:cone}
  \xymatrix@-.5em{
    A\ar[rr]^\varphi&&B\ar[dl]|-\circ\\
    &C_\varphi\ar[ul]
  }
\end{equation}
of a morphism \(\varphi\colon A\to B\) in \(\Cstarsep(X)\) is
defined as in~\cite{Meyer-Nest:BC} and is a diagram in
\(\KKcat(X)\).  The circled arrow from~\(B\) to~\(C_\varphi\)
means a \Star{}homomorphism \(\Sigma(B)\to C_\varphi\).  A
triangle in \(\KKcat(X)\) is called \emph{exact} if it is
isomorphic in \(\KKcat(X)\) to the mapping cone triangle of
some morphism in \(\Cstarsep(X)\).

As in~\cite{Meyer-Nest:BC}, there is an equivalent description
of the exact triangles using semi-split extensions in
\(\Cstarsep(X)\).  An extension
\begin{equation}
  \label{eq:extension}
  I\overset{i}\into E\overset{p}\prto Q
\end{equation}
gives rise to a commuting diagram
\[\xymatrix{
  \Sigma Q \ar@{=}[d]&
  I \ar[r]^{i}\ar[d]&
  E \ar[r]^{p}\ar@{=}[d]&
  Q \ar@{=}[d]\\
  \Sigma Q \ar[r]&
  C_p \ar[r]&
  E \ar[r]^{p}&
  Q.
}
\]

\begin{definition}
  \label{def:admissible}
  We call the extension \emph{admissible} if the map \(I\to
  C_p\) is invertible in \(\KKcat(X)\).
\end{definition}

The proof of the Excision Theorem~\ref{the:excision} shows that
this is the case if the extension is semi-split; but there are
more admissible extensions than semi-split extensions.  If the
extension is admissible, then there is a unique map \(\Sigma
Q\to I\) so that the top row becomes isomorphic to the bottom
row as a triangle in \(\KKcat(X)\).  Thus any admissible
extension in \(\Cstarsep(X)\) yields an exact triangle \(\Sigma
Q\to I\to E\to Q\), called \emph{extension triangle}.

Conversely, if \(\varphi\colon A\to B\) is a morphism in
\(\Cstarsep(X)\), then its mapping cone triangle is isomorphic
in \(\KKcat(X)\) to the extension triangle for the canonically
semi-split extension \(C_\varphi\into Z_\varphi\prto B\),
where~\(Z_\varphi\) denotes the mapping cylinder
of~\(\varphi\), which is homotopy equivalent to~\(A\).  The
above arguments work exactly as in the case of undecorated
Kasparov theory discussed in~\cite{Meyer-Nest:BC}.

As a result, a triangle in \(\KKcat(X)\) is isomorphic to a
mapping cone triangle of some morphism in \(\Cstarsep(X)\) if
and only if it is isomorphic to the extension triangle of some
semi-split extension in \(\Cstarsep(X)\).

\begin{proposition}
  \label{pro:KKX_triangulated}
  The category \(\KKcat(X)\) with the suspension automorphism
  and extension triangles specified above is a triangulated
  category.
\end{proposition}

\begin{proof}
  Most of the axioms amount to well-known properties of mapping
  cones and mapping cylinders, which are proven by translating
  corresponding arguments for the stable homotopy category of
  spaces, see~\cite{Meyer-Nest:BC}.

  The only axiom that requires a new argument in our case is
  (TR1), which asserts that any morphism in \(\KKcat(X)\) is
  part of some exact triangle.  The argument
  in~\cite{Meyer-Nest:BC} uses the description of Kasparov
  theory via the universal algebra \(qA\) by Joachim Cuntz.
  This approach can be made to work in \(\KKcat(X)\), but it is
  rather unflexible because the primitive ideal space of \(qA\)
  is hard to control.

  The following argument, which is inspired
  by~\cite{Bonkat:Thesis}, also applies to interesting
  subcategories of \(\KKcat(X)\) like the subcategory of
  nuclear \(\Cst\)\nb-algebras over~\(X\), which is studied
  in~\S\ref{sec:simple_fibres}.  Hence this is a triangulated
  category as well.

  Let \(f\in\KK_0(X;A,B)\).  We identify \(\KK_0(X;A,B) \cong
  \KK_1(X; A,\Sigma B)\).  Represent the image of~\(f\) in
  \(\KK_1(X; A, \Sigma B)\) by a cycle \((\varphi,\Hils,F)\).
  Adding a degenerate cycle, if necessary, we can achieve that
  the map \(\Phi\colon A\ni a\mapsto F^*\varphi(a)F\bmod
  \Comp(\Hils)\) is an injection from~\(A\) into the Calkin
  algebra \(\Bound(\Hils)\mathbin/ \Comp(\Hils)\) of~\(\Hils\)
  and that~\(\Hils\) is full, so that \(\Comp(\Hils)\) is
  \(\KK(X)\)\nb-equivalent to \(\Sigma B\).  The properties of
  a Kasparov cycle mean that~\(\Phi\) is the Busby invariant of
  a semi-split extension \(\Comp(\Hils) \into E\prto A\) of
  \(\Cst\)\nb-algebras over~\(X\).  The composition product of the
  map \(\Sigma A\to \Comp(\Hils)\) in the associated extension
  triangle and the canonical \(\KK(X)\)\nb-equivalence
  \(\Comp(\Hils)\simeq \Sigma B\) is the suspension of
  \(f\in\KK_0(X;A,B)\).  Hence we can embed~\(f\) in an exact
  triangle.
\end{proof}

\subsection{Adjointness relations}
\label{sec:adjointness}

\begin{proposition}
  \label{pro:adjointness_KK}
  Let~\(X\) be a topological space and let \(Y\in\Loclo(X)\).

  If \(Y\subseteq X\) is open, then we have natural
  isomorphisms
  \[
  \KK_*(X;i_Y^X(A),B) \cong \KK_*\bigl(Y;A,r_X^Y(B)\bigr)
  \]
  for all \(A\inOb\KKcat(Y)\), \(B\inOb\KKcat(X)\), that is,
  \(i_Y^X\) is left adjoint to \(r_X^Y\) as functors
  \(\KKcat(Y)\leftrightarrow\KKcat(X)\).

  If \(Y\subseteq X\) is closed, then we have natural
  isomorphisms
  \[
  \KK_*(Y;r_X^Y(A),B) \cong \KK_*\bigl(X;A,i_Y^X(B)\bigr)
  \]
  for all \(A\inOb\KKcat(X)\), \(B\inOb\KKcat(Y)\), that is,
  \(i_Y^X\) is right adjoint to \(r_X^Y\) as functors
  \(\KKcat(Y)\leftrightarrow\KKcat(X)\).
\end{proposition}

\begin{proof}
  Since both \(i_Y^X\) and \(r_X^Y\) descend to functors
  between \(\KKcat(X)\) and \(\KKcat(Y)\), this follows from
  the adjointness on the level of \(\Cstarcat(X)\) and
  \(\Cstarcat(Y)\) in Lemma~\ref{lem:adjointness}; an analogous
  assertion for induction and restriction functors for group
  actions on \(\Cst\)\nb-algebras is proven in
  \cite{Meyer-Nest:BC}*{\S3.2}.  The point of the argument is
  that an adjointness relation is equivalent to the existence
  of certain natural transformations called unit and counit of
  the adjunction, subject to some conditions
  (see~\cite{MacLane:Categories}).  These natural
  transformations already exist on the level of
  \Star{}homomorphisms, which induce morphisms in \(\KKcat(X)\)
  or \(\KKcat(Y)\).  The necessary relations for unit and
  counit of adjunction hold in \(\KKcat(\dotso)\) because they
  already hold in \(\Cstarcat(\dotso)\).  The unit and counit
  are natural in \(\KKcat(\dotso)\) and not just in
  \(\Cstarcat(\dotso)\) because of the uniqueness part of the
  universal property of \(\KKcat\).
\end{proof}

\begin{proposition}
  \label{pro:KK_ix}
  Let~\(X\) be a topological space and let \(x\in X\).  Then
  \[
  \KK_*\bigl(X;A,i_x(B)\bigr)
  \cong \KK_*\bigl(A\bigl(\cl{\{x\}}\bigr),B\bigr)
  \]
  for all \(A\inOb\Cstarsep(X)\), \(B\inOb\Cstarsep\).  That
  is, the functor \(i_x\colon \KKcat\to\KKcat(X)\) is right
  adjoint to the functor \(A\mapsto A\bigl(\cl{\{x\}}\bigr)\).
  Moreover,
  \[
  \KK_*(X;i_x(A),B)
  \cong \KK_*\biggl(A, \bigcap_{U\in\mathcal{U}_x} B(U)\biggr)
  \]
  for all \(A\inOb\Cstarsep\), \(B\inOb\Cstarsep(X)\),
  where~\(\mathcal{U}_x\) denotes the open neighbourhood filter
  of~\(x\) in~\(X\), That is, the functor \(i_x\colon
  \KKcat\to\KKcat(X)\) is left adjoint to the functor
  \(B\mapsto \bigcap_{U\in\mathcal{U}_x} B(U)\).  If~\(x\) has
  a minimal open neighbourhood~\(U_x\), then
  \[
  \KK_*(X;i_x(A),B) \cong \KK_*\bigl(A,B(U_x)\bigr),
  \]
\end{proposition}

\begin{proof}
  This follows from Lemma~\ref{lem:Hom_ix} in the same way as
  Proposition~\ref{pro:adjointness_KK}.  Notice that \(B\mapsto
  \bigcap_{U\in\mathcal{U}_x} B(U)\) commutes with
  \(\Cst\)\nb-stabilisation and maps (semi)-split extensions in
  \(\Cstarcat(X)\) again to (semi)-split extensions in
  \(\Cstarcat\); therefore, it descends to a functor
  \(\KKcat(X)\to\KKcat\).
\end{proof}

\section{The bootstrap class}
\label{sec:bootstrap}

Throughout this section, \(X\) denotes a \emph{finite} and
sober topological space.  Finiteness is crucial here.  First we
construct a canonical filtration on any \(\Cst\)\nb-algebra
over~\(X\).  We use this to study the analogue of the bootstrap
class in \(\KKcat(X)\).  Along the way, we also introduce the
larger category of local \(\Cst\)\nb-algebras over~\(X\).
Roughly speaking, locality means that all the canonical
\(\Cst\)\nb-algebra extensions that we get from
\(\Cst\)\nb-algebras over~\(X\) are admissible.  Objects in the
\(X\)\nb-equivariant bootstrap category have the additional
property that their fibres belong to the usual bootstrap
category.

\subsection{The canonical filtration}
\label{sec:canonical_filt}

We recursively construct a canonical filtration
\[
\emptyset=\Fil_0X \subset \Fil_1X\subset\dotsb\subset\Fil_\ell X=X
\]
of~\(X\) by open subsets~\(\Fil_jX\), such that the differences
\[
X_j \defeq \Fil_jX\setminus \Fil_{j-1}X
\]
are discrete for all \(j=1,\dotsc,\ell\).  In each step, we
let~\(X_j\) be the subset of all open points in
\(X\setminus\Fil_{j-1}X\) --~so that~\(X_j\) is discrete~-- and
put \(\Fil_jX=\Fil_{j-1}X\cup X_j\).  Equivalently, \(X_j\)
consists of all points of \(X\setminus\Fil_{j-1}X\) that are
maximal for the specialisation preorder~\(\prec\).  Since~\(X\)
is finite, \(X_j\) is non-empty unless \(\Fil_{j-1}X=X\), and
our recursion reaches~\(X\) after finitely many steps.

\begin{definition}
  \label{def:length}
  The \emph{length}~\(\ell\) of~\(X\) is the length of the
  longest chain \(x_1\prec x_2\prec\dotsb\prec x_\ell\)
  in~\(X\).
\end{definition}

We assume~\(X\) finite to ensure that the above filtration can
be constructed.  It is easy to extend our arguments to
\emph{Alexandrov spaces of finite length}; the only difference
is that the discrete spaces~\(X_j\) may be infinite in this
case, so that we need infinite direct sums in some places,
forcing us in Proposition~\ref{pro:local_over_X} to drop~(2)
and replace “triangulated” by “localising” in the last
sentence.  It should be possible to treat Alexandrov spaces of
infinite length in a similar way.  Since such techniques cannot
work for non-Alexandrov spaces, anyway, we do not pursue these
generalisations here.

\begin{definition}
  \label{def:PY}
  We shall use the functors
  \[
  P_Y \defeq i_Y^X\circ r_X^Y\colon \Cstarcat(X)\to\Cstarcat(X)
  \]
  for \(Y\in\Loclo(X)\).  Thus \((P_Y A)(Z) \cong A(Y\cap Z)\)
  for all \(Z\in\Loclo(X)\).
\end{definition}

If \(Y\in\Loclo(X)\), \(U\in\Open(Y)\), then we get an
extension
\begin{equation}
  \label{eq:PU_exact}
  P_U(A) \into P_Y(A) \prto P_{Y\setminus U}(A)
\end{equation}
in \(\Cstarcat(X)\) because of the extensions \(A(Z\cap U)
\into A(Z\cap Y) \prto A(Z\cap Y\setminus U)\) for all
\(Z\in\Loclo(X)\).

Let~\(A\) be a \(\Cst\)\nb-algebra over~\(X\).  We equip~\(A\) with
the canonical increasing filtration by the ideals
\[
\Fil_jA \defeq P_{\Fil_jX}(A),
\qquad
j=0,\dotsc,\ell,
\]
so that
\begin{equation}
  \label{eq:Filtration_over_X}
  \Fil_j A(Y) = A(Y\cap \Fil_jX) = A(Y)\cap A(\Fil_jX)
  \qquad
  \text{for all \(Y\in\Loclo(X)\).}
\end{equation}
Equation~\eqref{eq:PU_exact} shows that the subquotients of
this filtration are
\begin{equation}
  \label{eq:Filtration_over_X_subquotients}
  \Fil_j A\mathbin/ \Fil_{j-1} A
  \cong P_{\Fil_j X\setminus \Fil_{j-1} X}(A)
  = P_{X_j}(A)
  \cong \bigoplus_{x\in X_j} P_x(A)
  = \bigoplus_{x\in X_j} i_x\bigl(A(x)\bigr).
\end{equation}
Here \(i_x=i_x^X\) for \(x\in X\) denotes the extension functor
from the subset \(\{x\}\subseteq X\):
\[
i_x\colon \KKcat\congto \KKcat(\{x\}) \xrightarrow{i_x}
\KKcat(X),
\qquad
(i_x B)(Y) =
\begin{cases}
  B &\text{if \(x\in Y\),}\\
  0 &\text{if \(x\notin Y\).}
\end{cases}
\]

\begin{example}
  \label{exa:filtration_of_extension}
  Consider the space \(X=\{1,2\}\) with the non-discrete
  topology described in Example~\ref{exa:extensions}.  Here
  \[
  \Fil_0X=\emptyset,\quad
  \Fil_1X=\{1\},\quad
  \Fil_2X=\{1,2\}=X,\qquad
  X_1 =\{1\},\quad
  X_2 =\{2\}.
  \]
  The filtration \(\Fil_jA\) on a \(\Cst\)\nb-algebra over~\(X\)
  has one non-trivial layer \(\Fil_1A\) because
  \(\Fil_0A=\{0\}\) and \(\Fil_2A=A\).  Recall that
  \(\Cst\)\nb-algebras over~\(X\) correspond to extensions of
  \(\Cst\)\nb-algebras.  For a \(\Cst\)\nb-algebra extension \(I\into
  A\prto A/I\), the first filtration layer is simply the
  extension \(I\into I\prto 0\), so that the quotient
  \(A/\Fil_1A\) is the extension \(0\into A/I\prto A/I\).  Our
  filtration decomposes \(I\into A\prto A/I\) into an extension
  of \(\Cst\)\nb-algebra extensions as follows:
  \[
  (I\into I\prto 0) \into
  (I\into A\prto A/I) \prto
  (0\into A/I\prto A/I).
  \]
\end{example}

\begin{proposition}
  \label{pro:local_over_X}
  The following are equivalent for a separable
  \(\Cst\)\nb-algebra~\(A\) over~\(X\):
  \begin{enumerate}[label=\textup{(\arabic{*})}]
  \item The extensions \(\Fil_{j-1}A\into \Fil_jA\prto
    \Fil_jA\mathbin/ \Fil_{j-1}A\) in \(\Cstarsep(X)\) are
    admissible for \(j=1,\dotsc,\ell\).

  \item \(A\inOb\KKcat(X)\) belongs to the triangulated
    subcategory of \(\KKcat(X)\) generated by objects of the
    form \(i_x(B)\) with \(x\in X\), \(B\inOb\KKcat\).

  \item \(A\inOb\KKcat(X)\) belongs to the localising
    subcategory of \(\KKcat(X)\) generated by objects of the
    form \(i_x(B)\) with \(x\in X\), \(B\inOb\KKcat\).

  \item For any \(Y\in\Loclo(X)\), \(U\in\Open(Y)\), the
    extension
    \[
    P_U(A) \into P_Y(A) \prto P_{Y\setminus U}(A)
    \]
    in \(\Cstarsep(X)\) described above is admissible.
  \end{enumerate}
  Furthermore, if~\(A\) satisfies these conditions, then it
  already belongs to the triangulated subcategory of
  \(\KKcat(X)\) generated by \(i_x\bigl(A(x)\bigr)\) for \(x\in
  X\).
\end{proposition}

Recall that the \emph{localising subcategory} generated by a
family of objects in \(\KKcat(X)\) is the smallest subcategory
that contains the given objects and is triangulated and closed
under countable direct sums.

\begin{proof}
  (2)\(\Longrightarrow\)(3) and (4)\(\Longrightarrow\)(1) are
  trivial.  We will prove (1)\(\Longrightarrow\)(2) and
  (3)\(\Longrightarrow\)(4).
  \begin{description}
  \item[(1)\(\Longrightarrow\)(2)] Since the extensions
    \(\Fil_{j-1}A\into \Fil_jA \prto \Fil_jA\mathbin/
    \Fil_{j-1}A\) are admissible, they yield extension
    triangles in \(\KKcat(X)\).  Thus \(\Fil_jA\) belongs to
    the triangulated subcategory of \(\KKcat(X)\) generated by
    \(\Fil_{j-1}A\) and \(\Fil_jA\mathbin/ \Fil_{j-1}A\).  Since
    \({\Fil_0A=0}\), induction on~\(j\)
    and~\eqref{eq:Filtration_over_X_subquotients} show that
    \(\Fil_jA\) belongs to the triangulated subcategory
    generated by \(i_x A(x)\) with \(x\in\Fil_jX\).  Thus
    \(A=\Fil_\ell A\) belongs to the triangulated subcategory
    of \(\KKcat(X)\) generated by \(i_x\bigl(A(x)\bigr)\) for
    \(x\in X\).  This also yields the last statement in the
    proposition.

  \item[(3)\(\Longrightarrow\)(4)] It is clear that~(4) holds
    for objects of the form \(i_x(B)\) because at least one of
    the three objects \(P_U i_x(B)\), \(P_Y i_x(B)\), or
    \(P_{Y\setminus U}i_x(B)\) vanishes.  The property~(4) is
    inherited by (countable) direct sums, suspensions, and
    mapping cones.  To prove the latter, we use the definition
    of admissibility as an isomorphism statement in
    \(\KKcat(X)\) and the Five Lemma in triangulated
    categories.  Hence~(4) holds for all objects of the
    localising subcategory generated by \(i_x(B)\) for \(x\in
    X\), \(B\inOb\KKcat\).\qedhere
  \end{description}
\end{proof}

\begin{definition}
  \label{def:local_in_KK}
  Let \(\KKcat(X)_\loc\subseteq \KKcat(X)\) be the full
  subcategory of all objects that satisfy the equivalent
  conditions of Proposition~\ref{pro:local_over_X}.
\end{definition}

The functor \(f_*\colon \KKcat(X)\to\KKcat(Y)\) for a
continuous map \(f\colon X\to Y\) restricts to a functor
\(\KKcat(X)_\loc\to\KKcat(Y)_\loc\) because \(f_*\circ i_x^X=
i_x^Y\) and~\(f_*\) is an exact functor.  Similarly, the
restriction functor \(r_X^Y\colon \KKcat(X)\to\KKcat(Y)\) for a
locally closed subset \(Y\subseteq X\) maps \(\KKcat(X)_\loc\)
to \(\KKcat(Y)_\loc\) because it is exact and \(r_X^Y\circ
i_x^X\) is \(i_x^Y\) for \(x\in Y\) and~\(0\) otherwise.

\begin{proposition}
  \label{pro:local_invertible}
  Let~\(X\) be a finite topological space.  Let
  \(A,B\in\KKcat(X)_\loc\) and let \(f\in\KK_*(X;A,B)\).  If
  \(f(x)\in\KK_*\bigl(A(x),B(x)\bigr)\) is invertible for all
  \(x\in X\), then~\(f\) is invertible in \(\KKcat(X)\).  In
  particular, if \(A(x)\cong0\) in \(\KKcat\) for all \(x\in
  X\), then \(A\cong0\) in \(\KKcat(X)\).
\end{proposition}

\begin{proof}
  The second assertion follows immediately from the last
  sentence in Proposition~\ref{pro:local_over_X}.  It implies
  the first one by a well-known trick: embed~\(\alpha\) in an
  exact triangle by axiom (TR1) of a triangulated category, and
  use the long exact sequence to relate invertibility
  of~\(\alpha\) to the vanishing of its mapping cone.
\end{proof}

\begin{proposition}
  \label{pro:nuclear_local}
  Suppose that the \(\Cst\)\nb-algebra extensions
  \[
  A(U_x\setminus\{x\})\into A(U_x)\prto A(x)
  \]
  are semi-split for all \(x\in X\).  Then
  \(A\inOb\KKcat(X)_\loc\).  In particular, this applies if the
  underlying \(\Cst\)\nb-algebra of \(A\inOb\KKcat(X)\) is nuclear.
\end{proposition}

\begin{proof}
  We claim that the extensions in
  Proposition~\ref{pro:local_over_X}.(1) are semi-split as
  extensions of \(\Cst\)\nb-algebras over~\(X\), hence
  admissible in \(\KKcat(X)\).  For this, we need a completely
  positive section \(A(X_j)\to A(\Fil_j X)\) that is
  \(X\)\nb-equivariant, that is, restricts to maps \(A(X_j\cap
  V)\to A(\Fil_j X\cap V)\) for all \(V\in\Open(\Fil_j X)\).
  We take the sum of the completely positive sections for the
  extensions \(A(U_x\setminus\{x\})\into A(U_x)\prto A(x)\) for
  \(x\in X_j\).  This map has the required property because any
  open subset containing~\(x\) also contains~\(U_x\); it is
  irrelevant whether or not this section is contractive by the
  last sentence in Theorem~\ref{the:excision}.  If~\(A\) is
  nuclear, so are the ideals \(A(U_x)\) and their quotients
  \(A(x)\) for \(x\in X\).  Thus the above extensions have
  completely positive sections by the Choi--Effros Lifting
  Theorem (see~\cite{Choi-Effros:CP_lifting}).
\end{proof}

It is not clear whether the mere admissibility in~\(\Cstarsep\)
of the extensions
\[
A(U_x\setminus\{x\})\into A(U_x)\prto A(x)
\]
suffices to conclude that \(A\inOb\KKcat(X)_\loc\).  This
condition is certainly necessary.

\medskip

The constructions above yield spectral sequences as
in~\cite{Schochet:Top1}.  These may be useful for spaces of
length~\(1\), where they degenerate to a short exact sequence.
We only comment on this very briefly.

Let \(A\inOb\KKcat(X)_\loc\).  The admissible extensions
\(\Fil_{j-1} A\into \Fil_j A\prto \Fil_j A\mathbin/ \Fil_{j-1}
A\) for \(j=1,\dotsc,\ell\) produce exact triangles in
\(\KKcat(X)\).  A homological or cohomological functor such as
\(\KK(X;D,\blank)\) or \(\KK(X;\blank,D)\) maps these exact
triangles to a sequence of exact chain complexes.  These can be
arranged in an \emph{exact couple}, which generates a spectral
sequence (see~\cite{MacLane:Homology}).  This spectral sequence
could, in principle, be used to compute \(\KK_*(X;A,B)\) in
terms of
\[
\KK_*(X;\Fil_j A\mathbin/ \Fil_{j-1} A,B)
\cong \prod_{x\in X_j} \KK_*(X;i_x A(x),B)
\cong \prod_{x\in X_j} \KK_*\bigl(A(x),B(U_x)\bigr),
\]
where we have used Proposition~\ref{pro:KK_ix}.  These groups
comprise the \(E_1\)\nb-terms of the spectral sequence that we
get from our exact couple for the functor \(\KK(X;\blank,B)\).

For instance, consider again the situation of
Example~\ref{exa:filtration_of_extension}.  Let
\(I\triangleleft A\) and \(J\triangleleft B\) be
\(\Cst\)\nb-algebras over~\(X\), corresponding to \(\Cst\)\nb-algebra
extensions \(I\into A\prto A/I\) and \(J\into B\prto B/J\).
The above spectral sequence degenerates to a long exact
sequence
\[
\xymatrix{
  \KK_0(A/I,B) \ar[r]&
  \KK_0(X;I\triangleleft A,J\triangleleft B) \ar[r]&
  \KK_0(I,J) \ar[d]^{\delta}\\
  \KK_1(I,J) \ar[u]^{\delta}&
  \KK_1(X;I\triangleleft A,J\triangleleft B) \ar[l]&
  \KK_1(A/I,B). \ar[l]
}
\]
The boundary map is the diagonal map in the following commuting
diagram:
\[
\xymatrix{
  \KK_0(I,J) \ar[d] \ar[r] \ar[dr]^{\delta}&
  \KK_0(I,B) \ar[d]\\
  \KK_1(A/I,J) \ar[r]&
  \KK_1(A/I,B).
}
\]
We can rewrite the long exact sequence above as an extension:
\[
\coker\delta \into
\KK_*(X;I\triangleleft A,J\triangleleft B) \prto
\ker\delta.
\]
But we lack a description of \(\ker\delta\) and
\(\coker\delta\) as \(\Hom\)- and
\(\operatorname{Ext}\)-groups.  Therefore, the Universal
Coefficient Theorem of Alexander Bonkat~\cite{Bonkat:Thesis}
seems more attractive.

\subsection{The bootstrap class}
\label{sec:bootstrap_over_X}

The bootstrap class~\(\Bootstrap\) in \(\KKcat\) is the
localising subcategory generated by the single object~\(\C\),
that is, it is the smallest class of separable \(\Cst\)\nb-algebras
that contains~\(\C\) and is closed under \(\KK\)-equivalence,
countable direct sums, suspensions, and the formation of
mapping cones (see~\cite{Meyer-Nest:Filtrated_K}).

A localising subcategory of \(\KKcat(X)\) or \(\KKcat\) is
automatically closed under various other constructions, as
explained in~\cite{Meyer-Nest:BC}.  This includes admissible
extensions, admissible inductive limits (the appropriate notion
of admissibility is explained in~\cite{Meyer-Nest:BC}), and
crossed products by \(\Z\) and~\(\R\) and, more generally, by
actions of torsion-free amenable groups.

The latter result uses the reformulation of the (strong)
Baum--Connes property for such groups in~\cite{Meyer-Nest:BC}.
This reformulation asserts that~\(\C\) with the trivial
representation of an amenable group~\(G\) belongs to the
localising subcategory of \(\KKcat(G)\) generated by
\(\CONT_0(G)\).  Carrying this over to \(\KKcat(X)\), we
conclude that \(A\rtimes G\) for \(A\inOb\KKcat(X)\) belongs to
the localising subcategory of \(\KKcat(X)\) generated by
\(\bigl(A\otimes \CONT_0(G)\bigr)\rtimes G\), which is
Morita--Riefel equivalent to~\(A\).

The following definition provides an analogue
\(\Bootstrap(X)\subseteq \KKcat(X)\) of the bootstrap class
\(\Bootstrap\subseteq\KKcat\) for a finite topological
space~\(X\):

\begin{definition}
  \label{def:bootstrap}
  Let \(\Bootstrap(X)\) be the localising subcategory of
  \(\KKcat(X)\) that is generated by \(i_x(\C)\) for \(x\in
  X\).
\end{definition}

Notice that \(\{i_x(\C)\mid x\in X\}\) lists all possible ways
to turn~\(\C\) into a \(\Cst\)\nb-algebra over~\(X\).

\begin{proposition}
  \label{pro:bootstrap_over_X}
  Let~\(X\) be a finite topological space and let
  \(A\inOb\KKcat(X)\).  The following conditions are
  equivalent:
  \begin{enumerate}[label=\textup{(\arabic{*})}]
  \item \(A\inOb\Bootstrap(X)\);

  \item \(A\inOb\KKcat(X)_\loc\) and \(A(x)\inOb\Bootstrap\) for
    all \(x\in X\);

  \item the extensions \(\Fil_{j-1}A\into \Fil_jA\prto
    \Fil_jA\mathbin/ \Fil_{j-1}A\) are admissible for
    \(j=1,\dotsc,\ell\), and \(A(x)\inOb\Bootstrap\) for all
    \(x\in X\);

  \end{enumerate}
  In addition, in this case \(A(Y)\inOb\Bootstrap\) for all
  \(Y\in\Loclo(X)\).
\end{proposition}

\begin{proof}
  The equivalence of (2) and~(3) is already contained in
  Proposition~\ref{pro:local_over_X}.  Using the last sentence
  of Proposition~\ref{pro:local_over_X}, we also get the
  implication (3)\(\Longrightarrow\)(1) because~\(i_x\) is
  exact and commutes with direct sums.  The only asertion that
  is not yet contained in Proposition~\ref{pro:local_over_X} is
  that \(A\inOb\Bootstrap(X)\) implies \(A(Y)\inOb\Bootstrap\)
  for all \(Y\in\Loclo(X)\).  The reason is that the functor
  \(\KKcat(X)\to\KKcat\), \(A\mapsto A(Y)\), is exact,
  preserves countable direct sums, and maps the generators
  \(i_y(\C)\) for \(y\in X\) to either~\(0\) or~\(\C\) and
  hence into~\(\Bootstrap\).
\end{proof}

\begin{corollary}
  \label{cor:nuclear_bootstrap_over_X}
  If the underlying \(\Cst\)\nb-algebra of~\(A\) is nuclear, then
  \(A\inOb\Bootstrap(X)\) if and only if \(A(x)\inOb\Bootstrap\)
  for all \(x\in X\).
\end{corollary}

\begin{proof}
  Combine Propositions \ref{pro:nuclear_local}
  and~\ref{pro:bootstrap_over_X}.
\end{proof}

\begin{example}
  \label{exa:nuclear_bootstrap_over_X}
  View a separable nuclear \(\Cst\)\nb-algebra~\(A\) with only
  finitely many ideals as a \(\Cst\)\nb-algebra over \(\Prim(A)\).
  Example~\ref{exa:simple_subquotients} and
  Corollary~\ref{cor:nuclear_bootstrap_over_X} show that~\(A\)
  belongs to \(\Bootstrap(\Prim A)\) if and only if all its
  simple subquotients belong to the usual bootstrap class in
  \(\KKcat\).
\end{example}

\begin{proposition}
  \label{pro:bootstrap_invertible}
  Let~\(X\) be a finite topological space.  Let
  \(A,B\inOb\Bootstrap(X)\) and let \(f\in\KK_*(X;A,B)\).
  If~\(f\) induces invertible maps \(\K_*\bigl(A(x)\bigr)\to
  \K_*\bigl(B(x)\bigr)\) for all \(x\in X\), then~\(f\) is
  invertible in \(\KKcat(X)\).  In particular, if
  \(\K_*\bigl(A(x)\bigr)=0\) for all \(x\in X\), then
  \(A\cong0\) in \(\KKcat(X)\).
\end{proposition}

\begin{proof}
  As in the proof of Proposition~\ref{pro:local_invertible}, it
  suffices to show the second assertion.  Since
  \(A(x)\inOb\Bootstrap\) for all \(x\in X\), vanishing of
  \(\K_*\bigl(A(x)\bigr)\) implies vanishing of~\(A(x)\) in
  \(\KKcat\), so that Proposition~\ref{pro:local_invertible}
  yields the assertion.
\end{proof}

\subsection{Complementary subcategories}
\label{sec:complements}

It is often useful to replace a given object of \(\KKcat(X)\)
by one in the bootstrap class or \(\KKcat(X)_\loc\) that is as
close to the original as possible.  This is achieved by
localisation functors
\[
L_\Bootstrap\colon \KKcat(X)\to\Bootstrap(X),\qquad
L\colon \KKcat(X)\to\KKcat(X)_\loc
\]
that are right adjoint to the embeddings of these
subcategories.  That is, we want
\(\KK_*\bigl(X;A,L_\Bootstrap(B)\bigr) \cong \KK_*(X;A,B)\) for
all \(A\inOb\Bootstrap(X)\), \(B\in\KKcat(X)\) and similarly
for~\(L\).  These functors come with natural transformations
\(L_\Bootstrap \Rightarrow L\Rightarrow \ID\), and the defining
property is equivalent to \(L(A)_x\to A_x\) being a
\(\KK\)\nb-equivalence and \(\K_*(L_\Bootstrap(A)_x)\to
\K_*(A_x)\) being invertible for all \(x\in X\), respectively.

The functors \(L\) and~\(L_\Bootstrap\) exist because our two
subcategories belong to complementary pairs of localising
subcategories in the notation of~\cite{Meyer-Nest:BC}.  The
existence of this complementary pair is straightforward to
prove using the techniques of~\cite{Meyer:Homology_in_KK_II}.

\begin{definition}
  \label{def:complements}
  Let \(\KKcat(X)_\loc^\dashv\) be the class of all
  \(A\inOb\KKcat(X)\) for which \(A(x)\) is
  \(\KK\)\nb-equivalent to~\(0\) for all \(x\in X\).  Let
  \(\Bootstrap(X)^\dashv\) be the class of all
  \(A\inOb\KKcat(X)\) with \(\K_*\bigl(A(x)\bigr)=0\).
\end{definition}

\begin{theorem}
  \label{the:complementary}
  The pair of subcategories
  \((\Bootstrap(X),\Bootstrap(X)^\dashv)\) is complementary.
  So is the pair \((\KKcat(X)_\loc,\KKcat(X)_\loc^\dashv)\).
\end{theorem}

\begin{proof}
  We first prove the assertion for \(\KKcat(X)_\loc\).
  Consider the exact functor
  \[
  F\colon \KKcat(X)\to \prod_{x\in X} \KKcat,\qquad
  A\mapsto (A_x)_{x\in X}.
  \]
  Let~\(\Ideal\) be the kernel of~\(F\) on morphisms.
  Since~\(F\) is an exact functor that commutes with countable
  direct sums, \(\Ideal\) is a stable homological ideal that is
  compatible with direct sums
  (see~\cites{Meyer-Nest:Homology_in_KK,
    Meyer:Homology_in_KK_II}).  The kernel of~\(F\) on objects
  is exactly \(\KKcat(X)^\dashv\).  Proposition~\ref{pro:KK_ix}
  shows that the functor~\(F\) has a left adjoint, namely, the
  functor
  \[
  F^\lad\bigl((A_x)_{x\in X}\bigr) \defeq \bigoplus_{x\in X}
  i_x(A_x).
  \]
  Therefore, the ideal~\(\Ideal\) has enough projective objects
  by \cite{Meyer-Nest:Homology_in_KK}*{Proposition 3.37};
  furthermore, the projective objects are retracts of direct
  sums of objects of the form \(i_x(A_x)\).  Hence the
  localising subcategory generated by the
  \(\Ideal\)\nb-projective objects is \(\KKcat(X)_\loc\) by
  Proposition~\ref{pro:local_over_X}.  Finally,
  \cite{Meyer:Homology_in_KK_II}*{Theorem 4.6} shows that the
  pair of subcategories
  \((\KKcat(X)_\loc,\KKcat(X)_\loc^\dashv)\) is complementary.

  The argument for the bootstrap category is almost literally
  the same, but using the stable homological functor
  \(\K_*\circ F\colon \KKcat(X)\to \prod_{x\in X}
  \Ab^{\Z/2}_\textup{c}\) instead of~\(F\), where
  \(\Ab^{\Z/2}_\textup{c}\) denotes the category of countable
  \(\Z/2\)-graded Abelian groups.  The adjoint of \(\K_*\circ
  F\) is defined on families of countable free Abelian groups,
  which is enough to conclude that \(\ker (\K_*\circ F)\) has
  enough projective objects.  This time, the projective objects
  generate the category \(\Bootstrap(X)\), and the kernel
  of~\(F\) on objects is \(\Bootstrap(X)^\dashv\).  Hence
  \cite{Meyer:Homology_in_KK_II}*{Theorem 4.6} shows that the
  pair of subcategories
  \((\Bootstrap(X),\Bootstrap(X)^\dashv)\) is complementary.
\end{proof}

\begin{lemma}
  \label{lem:bootstrap_ortho_finite}
  The following are equivalent for \(A\inOb\KKcat(X)\):
  \begin{enumerate}[label=\textup{(\arabic{*})}]
  \item \(\K_*\bigl(A(x)\bigr)=0\) for all \(x\in X\), that is,
    \(A\inOb\Bootstrap(X)^\dashv\);

  \item \(\K_*\bigl(A(Y)\bigr)=0\) for all \(Y\in\Loclo(X)\);

  \item \(\K_*\bigl(A(U)\bigr)=0\) for all \(U\in\Open(X)\).
  \end{enumerate}
\end{lemma}

\begin{proof}
  It is clear that (2) implies both (1) and~(3).  Conversely,
  (3) implies~(2): write \(Y\in\Loclo(X)\) as \(U\setminus V\)
  with \(U,V\in\Open(X)\), \(V\subseteq U\), and use the
  \(\K\)\nb-theory long exact sequence for the extension
  \(A(U)\into A(V)\prto A(Y)\).  It remains to check that~(1)
  implies~(2).

  We prove by induction on~\(j\) that~(1) implies
  \(\K_*\bigl(A(Y)\bigr)=0\) for all \(Y\in\Loclo(\Fil_jX)\).
  This is trivial for \(j=0\).  If \(Y\subseteq\Fil_{j+1}X\),
  then \(\K_*\bigl(A(Y\cap\Fil_jX)\bigr)=0\) by the induction
  assumption.  The \(\K\)\nb-theory long exact sequence for the
  extension
  \[
  A(Y\cap\Fil_jX) \into A(Y)\prto
  \bigoplus_{x\in X_{j+1}\cap Y} A(x)
  \]
  yields \(\K_*\bigl(A(Y)\bigr)=0\) as claimed.
\end{proof}

We can also apply the machinery of
\cites{Meyer-Nest:Homology_in_KK, Meyer:Homology_in_KK_II} to
the ideal~\(\Ideal\) to generate a spectral sequence that
computes \(\KK_*(X;A,B)\).  This spectral sequence is more
useful than the one from the canonical filtration because its
second page involves derived functors.  But this spectral
sequence rarely degenerates to an exact sequence.

\subsection{A definition for infinite spaces}
\label{sec:def_bootstrap_infinite}

The ideas in~\S\ref{sec:complements} suggest a definition of
the bootstrap class for infinite spaces.

\begin{definition}
  \label{def:bootstrap_ortho_infinite}
  Let~\(X\) be a topological space.  Let
  \(\Bootstrap(X)^\dashv\subseteq\KKcat(X)\) consist of all
  separable \(\Cst\)\nb-algebras over~\(X\) with
  \(\K_*\bigl(A(U)\bigr)=0\) for all \(U\in\Open(X)\).
\end{definition}

Lemma~\ref{lem:bootstrap_ortho_finite} shows that this agrees
with our previous definition for finite~\(X\).  Furthermore,
the same argument as in the proof of
Lemma~\ref{lem:bootstrap_ortho_finite} yields
\(A\inOb\Bootstrap(X)^\dashv\) if and only if
\(\K_*\bigl(A(Y)\bigr)=0\) for all \(Y\in\Loclo(X)\).  The
first condition in Lemma~\ref{lem:bootstrap_ortho_finite} has
no analogue because of Example~\ref{exa:bad_interval}.

It is clear from the definition that \(\Bootstrap(X)^\dashv\)
is a localising subcategory of \(\KKcat(X)\).

\begin{definition}
  \label{def:bootstrap_infinite}
  Let~\(X\) be a topological space.  We let \(\Bootstrap(X)\)
  be the localisation of \(\KKcat(X)\) at
  \(\Bootstrap(X)^\dashv\).
\end{definition}

For finite~\(X\), we have seen that \(\Bootstrap(X)^\dashv\) is
part of a complementary pair of localising subcategories, with
partner \(\Bootstrap(X)\).  This shows that the localisation of
\(\KKcat(X)\) at \(\Bootstrap(X)^\dashv\) is canonically
equivalent to \(\Bootstrap(X)\).  For infinite~\(X\), it is
unclear whether \(\Bootstrap(X)^\dashv\) is part of a
complementary pair.  If it is, the partner must be
\[
\Bootstrapc \defeq \{A\inOb\KKcat(X)\mid
\text{\(\KK_*(X;A,B)=0\) for all
  \(B\inOb\Bootstrap(X)^\dashv\)}\}.
\]

Since \(B\inOb\Bootstrap(X)^\dashv\) implies nothing about the
\(\K\)\nb-theory of \(\bigcap_{U\in\mathcal{U}_x} B(U)\), in
general, Proposition~\ref{pro:KK_ix} shows that \(i_x\C\) does
not belong to~\(\Bootstrapc\) in general.

If~\(X\) is Hausdorff, then
\(\CONT_0(U)\inOb\Bootstrap(X)^\dashv\) for all
\(U\in\Open(X)\).  Nevertheless, it is not clear whether
\((\Bootstrapc,\Bootstrap(X)^\dashv)\) is complementary.

\section{Making the fibres simple}
\label{sec:simple_fibres}

\begin{definition}
  \label{def:tight}
  A \(\Cst\)\nb-algebra \((A,\psi)\) over~\(X\) is called
  \emph{tight} if \(\psi\colon \Prim(A)\to X\) is a
  homeomorphism.
\end{definition}

Tightness implies that the fibres \(A_x=A(x)\) for \(x\in X\)
are simple \(\Cst\)\nb-algebras.  But the converse does not
hold: the fibres are simple if and only if the map \(\psi\colon
\Prim(A)\to X\) is bijective.

To equip \(\KKcat(X)\) with a triangulated category structure,
we must drop the tightness assumption because it is usually
destroyed when we construct cylinders, mapping cones, or
extensions of \(\Cst\)\nb-algebras over~\(X\).  Nevertheless, we
show below that we may reinstall tightness by passing to a
\(\KK(X)\)-equivalent object, at least in the nuclear case.

The special case where the space~\(X\) in question has only one
point is already known:

\begin{theorem}[\cite{Rordam-Stormer:Classification_Entropy}*{Proposition 8.4.5}]
  \label{the:exhaust}
  Any separable nuclear \(\Cst\)\nb-algebra is
  \(\KK\)-equivalent to a \(\Cst\)\nb-algebra that is separable,
  nuclear, purely infinite, \(\Cst\)\nb-stable and simple.
\end{theorem}

Stability is not part of the assertion
in~\cite{Rordam-Stormer:Classification_Entropy}, but can be
achieved by tensoring with the compact operators, without
destroying the other properties.  The main difficulty is to
achieve simplicity.  We are going to generalise this theorem as
follows:

\begin{theorem}
  \label{the:exhaust_X}
  Let~\(X\) be a finite topological space.  Any separable
  nuclear \(\Cst\)\nb-algebra over~\(X\) is
  \(\KK(X)\)\nb-equivalent to a \(\Cst\)\nb-algebra over~\(X\)
  that is tight, separable, nuclear, purely infinite, and
  \(\Cst\)\nb-stable.
\end{theorem}

For the zero \(\Cst\)\nb-algebra, viewed as a
\(\Cst\)\nb-algebra over~\(X\), this reproves the known
statement that there is a separable, nuclear, purely infinite,
and stable \(\Cst\)\nb-algebra with spectrum~\(X\) for any
finite topological space~\(X\).

\begin{proof}
  Since~\(A\) is separable and nuclear, so are the
  subquotients~\(A_x\).  Hence Theorem~\ref{the:exhaust}
  provides simple, separable, nuclear, stable, purely infinite
  \(\Cst\)\nb-algebras \(B_x\) and \(\KK\)\nb-equivalences
  \(f_x\in\KK_0(A_x,B_x)\) for all \(x\in X\).

  We use the canonical filtration \(\Fil_j X\) of~\(X\) and the
  resulting filtration \(\Fil_j A\) introduced
  in~\S\ref{sec:canonical_filt},
  see~\eqref{eq:Filtration_over_X}.  The subquotients
  \[
  A_j^0 \defeq \Fil_j A\mathbin/ \Fil_{j-1} A
  \]
  of the filtration are described
  in~\eqref{eq:Filtration_over_X_subquotients} in terms of the
  subquotients~\(A_x\) for \(x\in X\).

  We will recursively construct a sequence~\(B_j\) of
  \(\Cst\)\nb-algebras over~\(X\) that are supported on~\(\Fil_j
  X\) and \(\KK(X)\)\nb-equivalent to~\(\Fil_j A\) for
  \(j=0,\dotsc,\ell\), such that \(\Fil_j B_k=B_j\) for \(k\ge
  j\) and each~\(B_j\) is tight over~\(\Fil_j X\), separable,
  nuclear, purely infinite, and stable.  The last
  object~\(B_\ell\) in this series is \(\KK(X)\)\nb-equivalent
  to \(\Fil_\ell A = A\) and has all the required properties.
  Since \(\Fil_0 X=\emptyset\), the recursion must begin with
  \(B_0=A_0=\{0\}\).  We assume that~\(B_j\) has been
  constructed.  Let
  \[
  B_{j+1}^0 \defeq \bigoplus_{x\in X_{j+1}} i_x(B_x).
  \]
  We will construct~\(B_{j+1}\) as an extension of~\(B_j\)
  by~\(B_{j+1}^0\).  This ensures that the fibres of~\(B_j\)
  are~\(B_x\) for \(x\in\Fil_j X\) and~\(0\) for \(x\in
  X\setminus \Fil_j X\).

  First we construct, for each \(x\in X_{j+1}\), a suitable
  extension of \(B_x\) by~\(B_j\).  Let \(U_x\subseteq
  \Fil_{j+1} X\) be the minimal open subset containing~\(x\)
  and let \(U_x'\defeq U_x\setminus\{x\}\).  Since~\(X_{j+1}\)
  is discrete, \(U_x'\) is an open subset of~\(\Fil_jX\).  The
  extension
  \[
  A(U_x') \into A(U_x) \prto A_x
  \]
  is semi-split and thus provides a class~\(\delta_x^A\) in
  \(\KK_1\bigl(A_x,A(U_x')\bigr)\) because~\(A_x\) is nuclear.
  Since \(B_x\simeq A_x\), \(\Fil_j A\simeq B_j\) and \(\Fil_j
  A(U_x') = A(U_x')\), we can transform this class to
  \(\delta_x^B\) in \(\KK_1\bigl(B_x,B_j(U_x')\bigr)\).

  We abbreviate \(B_{jx} \defeq B_j(U_x')\) to simplify our
  notation.  Represent \(\delta_x^B\) by an odd Kasparov cycle
  \((\Hils,\varphi,F)\), where~\(\Hils\) is a Hilbert
  \(B_{jx}\)\nb-module, \(\varphi\colon B_x \to \Bound(\Hils)\)
  is a \Star{}homomorphism, and \(F\in\Bound(\Hils)\) satisfies
  \(F^2=1\), \(F=F^*\), and \([F,\varphi(b)]\in \Comp(\Hils)\)
  for all \(b\in B_x\).  Now we apply the familiar
  correspondence between odd \(\KK\)\nb-elements and
  \(\Cst\)\nb-algebra extensions.  Let \(P\defeq
  \frac{1}{2}(1+F)\), then
  \[
  \psi\colon B_x \to \Bound(\Hils)\mathbin/ \Comp(\Hils),
  \qquad b\mapsto P\varphi(x)P
  \]
  is a \Star{}homomorphism and hence the Busby invariant of an
  extension of~\(B_x\) by \(\Comp(\Hils)\).  After adding a
  sufficiently big split extension, that is, a
  \Star{}homomorphism \(\psi_0\colon B_x\to \Bound(\Hils')\),
  the map \(\psi\colon B_x \to \Bound(\Hils)\mathbin/
  \Comp(\Hils)\) becomes injective and the ideal in
  \(\Comp(\Hils)\) generated by
  \(\Comp(\Hils)\psi(B_x)\Comp(\Hils)\) is all of
  \(\Comp(\Hils)\).  We assume these two extra properties from
  now on.

  We also add to~\(\psi\) the trivial extension \(B_j\into
  B_j\oplus B_x\prto B_x\), whose Busby invariant is the zero
  map.  This produces an extension of~\(B_x\) by
  \(\Comp(B_j\oplus\Hils) \cong B_j\); the last isomorphism
  holds because~\(B_j\) is stable, so that
  \(B_j\oplus\Hils'\cong B_j\) for any Hilbert
  \(B_j\)\nb-module~\(\Hils'\).  Since~\(\psi\) is injective,
  the extension we get is of the form \(B_j\into E_{jx}\prto
  B_x\).  This extension is still semi-split, and its class in
  \(\KK_1(B_x,B_j)\) is the composite of \(\delta_x^B\) with
  the embedding \(B_{jx}\to B_j\).  Our careful construction
  ensures that the ideal in~\(B_j\) generated by
  \(B_j\psi(B_x)B_j\) is equal to \(B(U_x')\).

  Now we combine these extensions for all \(x\in X\) by taking
  their external direct sum.  This is an extension of
  \(\bigoplus_{x\in X_{j+1}} B_x = B_{j+1}^0\) by the
  \(\Cst\)\nb-algebra of compact operators on the Hilbert
  \(B_j\)\nb-module \(\bigoplus_{x\in X_{j+1}} B_j \cong B_j\),
  where we used the stability of~\(B_j\) once more.  Thus we
  obtain an extension \(B_j \into B_{j+1} \prto B_{j+1}^0\).
  We claim that the primitive ideal space of~\(B_{j+1}\)
  identifies naturally with~\(\Fil_{j+1} X\).

  The extension \(B_j\into B_{j+1}\prto B_{j+1}^0\) decomposes
  \(\Prim(B_{j+1})\) into an open subset homeomorphic to
  \(\Prim(B_j) \cong \Fil_j X\) and a closed subset
  homeomorphic to the discrete set \(\Prim(B_{j+1}^0) =
  X_{j+1}\).  This provides a canonical bijection between
  \(\Prim(B_{j+1})\) and \(\Fil_{j+1} X\).  We must check that
  it is a homeomorphism.

  First let \(U\subseteq \Fil_{j+1}X\) be open in
  \(\Fil_{j+1}X\).  Then \(U\cap\Fil_jX\) is open and
  contains~\(U_x'\) for each \(x\in U \cap X_{j+1}\).  Our
  construction ensures that \(\psi(B_x)\subseteq B_{j+1}\)
  multiplies~\(B_j\) into \(B_{jx}\subseteq B_j(U\cap \Fil_j
  X)\).  Hence \(B_j(U\cap\Fil_jX) + \sum_{x\in U\cap X_{j+1}}
  \psi(B_x)\) is an ideal in~\(B_{j+1}\).  This shows
  that~\(U\) is open in \(\Prim(B_{j+1})\).

  Now let \(U\subseteq \Fil_{j+1}X\) be open in
  \(\Prim(B_{j+1})\).  Then \(U\cap\Fil_j\) must be open in
  \(\Fil_jX\cong \Prim(B_j)\).  Furthermore, if \(x\in U\cap
  X_{j+1}\), then the subset of \(\Prim(B_j)\) corresponding to
  the ideal in~\(B_j\) generated by \(B_j\psi(B_x)B_j\) is
  contained in~\(U\).  But our construction ensures that this
  subset is precisely~\(U_x'\).  Hence
  \[
  U= (U\cap\Fil_j X) \cup \bigcup_{x\in U\cap X_{j+1}} U_x,
  \]
  proving that~\(U\) is open in the topology
  of~\(\Fil_{j+1}X\).  This establishes that our canonical map
  between \(\Prim(B_{j+1})\) and \(\Fil_{j+1}X\) is a
  homeomorphism.  Thus we may view~\(B_{j+1}\) as a
  \(\Cst\)\nb-algebra over~\(X\) supported in \(\Fil_{j+1}X\).
  It is clear from our construction that \(B_j\into
  B_{j+1}\prto B_{j+1}^0\) is an extension of
  \(\Cst\)\nb-algebras over~\(X\).  Here we view \(B_{j+1}^0\)
  as a \(\Cst\)\nb-algebra over~\(X\) in the obvious way, so
  that~\(B_x\) is its fibre over~\(x\) for \(x\in X_{j+1}\).

  There is no reason to expect \(B_{j+1}\) to be stable or
  purely infinite.  But this is easily repaired by tensoring
  with \(\Comp\otimes\mathcal{O}_\infty\).  This does not
  change \(B_j\) and~\(B_{j+1}^0\), up to isomorphism, because
  these are already stable and purely infinite, and it has no
  effect on the primitive ideal space, nuclearity or
  separability.  Thus we may achieve that~\(B_{j+1}\) is stable
  and purely infinite.

  By assumption, there is a \(\KK(X)\)\nb-equivalence
  \(f_j\in\KK_0(X;\Fil_j A,B_j)\).  Furthermore, our
  construction of~\(B_{j+1}^0\) ensures a
  \(\KK(X)\)\nb-equivalence \(f_{j+1}^0\) between \(A_{j+1}^0\)
  and~\(B_{j+1}^0\).  Due to the nuclearity of~\(A\), the
  arguments in~\S\ref{sec:canonical_filt} show that
  \[
  \Fil_j A\into \Fil_{j+1} A \prto A_{j+1}^0
  \]
  is a semi-split extension of \(\Cst\)\nb-algebras over~\(X\)
  and hence provides an exact triangle in \(\KKcat(X)\).  The
  same argument provides an extension triangle for the
  extension \(B_j\into B_{j+1}\prto B_{j+1}^0\).  Let
  \(\delta_j^A\) and~\(\delta_j^B\) be the classes in
  \(\KK_1(X;A_{j+1}^0, \Fil_j A)\) and \(\KK_1(B_{j+1}^0,B_j)\)
  associated to these extension; they appear in the exact
  triangles described above.

  Both classes \(\delta_j^A\) and \(\delta_j^B\) are,
  essentially, the sum of the classes \(\delta_x^A\) and
  \(\delta_x^B\) for \(x\in X_{j+1}\), respectively.  More
  precisely, we have to compose each \(\delta_x^A\) with the
  embedding \(A(U_x')\to \Fil_j A\).  Hence the solid square in
  the diagram
  \[
  \xymatrix{
    \Sigma A_{j+1}^0 \ar[r]^{\delta_j^A}\ar[d]_{\Sigma f_{j+1}^0}^{\cong}&
    \Fil_j A \ar[r]\ar[d]_{f_j}^{\cong}&
    \Fil_{j+1}A \ar[r]\ar@{..>}[d]_{f_{j+1}}^{\cong}&
    A_{j+1}^0 \ar[d]_{f_{j+1}^0}^{\cong}\\
    \Sigma B_{j+1}^0 \ar[r]^{\delta_j^B}&
    B_j \ar[r]&B_{j+1}\ar[r]&B_{j+1}^0
  }
  \]
  commutes.  By an axiom of triangulated categories, we can
  find the dotted arrow making the whole diagram commute.  The
  Five Lemma for triangulated categories asserts that this
  arrow is invertible because \(f_j\) and~\(f_{j+1}^0\) are.
  This shows that~\(B_{j+1}\) has all required properties and
  completes the induction step.
\end{proof}

\begin{theorem}
  \label{the:tighten_local}
  Let~\(X\) be a finite topological space and let~\(A\) be a
  separable \(\Cst\)\nb-algebra over~\(X\).  The following are
  equivalent:
  \begin{itemize}
  \item \(A\inOb\KKcat(X)_\loc\) and~\(A_x\) is
    \(\KK\)\nb-equivalent to a nuclear \(\Cst\)\nb-algebra for
    each \(x\in X\);
  \item \(A\) is \(\KK(X)\)\nb-equivalent to a
    \(\Cst\)\nb-algebra over~\(X\) that is tight, separable,
    nuclear, purely infinite, and \(\Cst\)\nb-stable.
  \end{itemize}
\end{theorem}

\begin{proof}
  The proof of Theorem~\ref{the:exhaust_X} still works under
  the weaker assumption that \(A\inOb\KKcat(X)_\loc\)
  and~\(A_x\) is \(\KK\)\nb-equivalent to a nuclear
  \(\Cst\)\nb-algebra for each \(x\in X\).  The converse
  implication is trivial.
\end{proof}

\begin{corollary}
  \label{cor:tighten_local}
  Let~\(X\) be a finite topological space and let~\(A\) be a
  separable \(\Cst\)\nb-algebra over~\(X\).  The following are
  equivalent:
  \begin{itemize}
  \item \(A\inOb\Bootstrap(X)\);

  \item \(A\) is \(\KK(X)\)\nb-equivalent to a
    \(\Cst\)\nb-algebra over~\(X\) that is tight, separable,
    nuclear, purely infinite, \(\Cst\)\nb-stable, and has
    fibres~\(A_x\) in the bootstrap class~\(\Bootstrap\).
  \end{itemize}
\end{corollary}

\begin{proof}
  Combine Theorem~\ref{the:tighten_local} and
  Proposition~\ref{pro:bootstrap_over_X}.
\end{proof}

By a deep classification result by Eberhard Kirchberg
(see~\cite{Kirchberg:Michael}), two tight, separable, nuclear,
purely infinite, stable \(\Cst\)\nb-algebras over~\(X\) are
\(\KK(X)\)-equivalent if and only if they are isomorphic as
\(\Cst\)\nb-algebras over~\(X\).  Therefore, the
representatives found in Theorems \ref{the:exhaust_X}
and~\ref{the:tighten_local} are unique up to
\(X\)\nb-equivariant \Star{}isomorphism.

Let \(\KKcat(X)_\textup{nuc}\) be the subcategory of
\(\KKcat(X)\) whose objects are the separable nuclear
\(\Cst\)\nb-algebras over~\(X\).  This is a triangulated
category as well because the basic constructions like
suspensions, mapping cones, and extensions never leave this
subcategory.  The subcategory of \(\KKcat(X)\) whose objects
are the tight, separable, nuclear, purely infinite, stable
\(\Cst\)\nb-algebras over~\(X\) is equivalent to
\(\KKcat(X)_\textup{nuc}\) by Theorem~\ref{the:exhaust_X} and
hence inherits a triangulated category structure.  It has the
remarkable feature that isomorphisms in this triangulated
category lift to \(X\)\nb-equivariant \Star{}isomorphisms.

Recall that a \(\Cst\)\nb-algebra belongs to \(\Bootstrap\) if
and only if it is \(\KK\)\nb-equivalent to a \emph{commutative}
\(\Cst\)\nb-algebra.  This probably remains the case at least
for finite spaces~\(X\), but the authors do not know how to
prove this.  For infinite spaces, it is even less clear whether
\(\Bootstrap(X)\) is equivalent to the \(\KK(X)\)-category of
commutative \(\Cst\)\nb-algebras over~\(X\).  We only have the
following characterisation:

\begin{theorem}
  \label{the:bootstrap_type_one}
  A separable \(\Cst\)\nb-algebra over~\(X\) belongs to the
  bootstrap class \(\Bootstrap(X)\) if and only if it is
  \(\KK(X)\)-equivalent to a \(\Cst\)\nb-stable, separable
  \(\Cst\)\nb-algebra over~\(X\) of type~\(I\).
\end{theorem}

\begin{proof}
  Follow the proof of Theorem~\ref{the:exhaust_X}, but using
  stabilisations of commutative \(\Cst\)\nb-algebras~\(B_x\)
  instead of nuclear purely infinite ones.  The proof shows
  that we can also achieve that the fibres~\(B_x\) are all of
  the form \(\CONT_0(Y_x)\otimes\Comp\) for second countable
  locally compact spaces~\(Y_x\).
\end{proof}

\section{Outlook}
\label{sec:outlook}

We have defined a bootstrap class
\(\Bootstrap(X)\subseteq\KKcat(X)\) over a finite topological
space~\(X\), which is the domain on which we should expect a
Universal Coefficient Theorem to compute \(\KK_*(X;A,B)\).  We
have seen that any object of the bootstrap class is
\(\KK(X)\)-equivalent to a tight, purely infinite, stable,
nuclear, separable \(\Cst\)\nb-algebra over~\(X\), for which
Kirchberg's classification results apply.

There are several spectral sequences that compute
\(\KK_*(X;A,B)\), but applications to the classification
programme require a short exact sequence.  For some finite
topological spaces, such a short exact sequence is constructed
in~\cite{Meyer-Nest:Filtrated_K} based on filtrated
\(\K\)\nb-theory, so that filtrated \(\K\)\nb-theory is a
complete invariant.  This invariant comprises the
\(\K\)\nb-theory \(\K_*\bigl(A(Y)\bigr)\) of all locally closed
subsets~\(Y\) of~\(X\) together with the action of all natural
transformations between them.  This is a consequence of a
Universal Coefficient Theorem in this case.  It is also shown
in~\cite{Meyer-Nest:Filtrated_K} that there are finite
topological spaces for which filtrated \(\K\)\nb-theory is not
yet a complete invariant.  At the moment, it is unclear whether
there is a general, tractable complete invariant for objects of
\(\Bootstrap(X)\).

Another issue is to treat infinite topological spaces.  A
promising approach is to approximate infinite spaces by finite
non-Hausdorff spaces associated to open coverings of the space
in question.  In good cases, there should be a
\(\varprojlim\nolimits^1\)-sequence that relates
\(\KK_*(X;A,B)\) to Kasparov groups over such finite
approximations to~\(X\), reducing computations from the
infinite to the finite case.  Such an exact sequence may be
considerably easier for \(X\)\nb-equivariant \(E\)\nb-theory,
where we do not have to worry about completely positive
sections.

\end{document}